\documentclass{article}
\usepackage[utf8]{inputenc}
\usepackage[a4paper, top=2.25cm, bottom=2.25cm, left=2.25cm, right=2.25cm]{geometry}
\usepackage{multicol}
\usepackage{amsmath}
\usepackage{siunitx}
\usepackage{booktabs}
\usepackage{graphicx}
\usepackage{tikz}
\definecolor{warmyellow}{RGB}{253, 243, 207}
\definecolor{sagegreen}{RGB}{209, 228, 218}
\usetikzlibrary{decorations.pathreplacing}
\usepackage{float}
\usepackage{dblfloatfix}
\usepackage{caption}
\usepackage{xcolor}
\definecolor{myorange}{RGB}{230, 130, 0}
\definecolor{darkgreen}{RGB}{60, 160, 60}
\usepackage{csquotes}
\usepackage[backend=bibtex,style=authoryear,citestyle=authoryear]{biblatex}
\usepackage{xurl}
\usepackage{hyperref}
\hypersetup{
breaklinks=true,
colorlinks=true,
linkcolor=black,
citecolor=black,
urlcolor=black
}

\usepackage{orcidlink}
\raggedcolumns
\begin{document}
\title{\vspace{-1cm}
The Ishango Bone: Evidence for Intentional Arithmetic Design \\ in the Late Stone Age?}
\author{Jenny Baur\,\orcidlink{0009-0009-8855-211X}\\
{\small{jenny.baur@outlook.com}}}
\date{}
\maketitle
\setlength{\parindent}{0pt}
\setlength{\parskip}{0.5\baselineskip}

\begin{center}
\begin{minipage}{0.9\linewidth}
\setlength{\parindent}{0pt} 
\setlength{\parskip}{0.5\baselineskip} 
The Ishango Bone is a prehistoric artefact dated to approximately 25,000--20,000 years ago, discovered near the Semliki River in what is now the Democratic Republic of Congo \parencite{Crevecoeur2016, ICOMOS2010}. It displays sixteen groups of notches organised into three columns, yet its precise function remains debated. This study identifies arithmetic patterns across all three columns through a sequence of exploratory observations and adjustments. Columns G and D together contain every odd number between 9 and 21 except 15, which is the arithmetic mean of the numbers represented; each column sums to 60 and contains two pairs that each sum to 30. Exchanging 11 and 21 in Column D aligns its grouping with that of Column G and yields further patterns. An analogous exchange of the first two pairs in Column M, together with the provisional reading $(9,4)$ for the damaged pair and a total-preserving
transfer of one unit from 5 to 3, brings the products of the M pairs into exact agreement with the G--D row sums. A common grouping rule can then be described across all three columns. Two conditional quantitative analyses, a Monte Carlo randomisation and a complete enumeration of 64 permitted readings and adjustments, characterise this structure within specified models. The adjusted arrangement offers reasonable grounds for considering intentional arithmetic organisation a possibility. This study aims to motivate further research on this bone and on other artefacts.
\end{minipage}
\end{center}

\vspace{0.5cm}

\section{Introduction}
\label{sec:intro}

The Ishango Bone is a Late Stone Age artefact, unearthed in 1950 by Belgian geologist Jean de Heinzelin during excavations along the Semliki River near Lake Edward in what is now the Democratic Republic of Congo \parencite{deHeinzelin1957}. Approximately 10 centimetres in length, with a piece of quartz affixed to one end, the bone is engraved with sixteen groups of notches organised into three distinct columns, an arrangement that has intrigued researchers across disciplines for decades and is often cited as one of the earliest known mathematical artefacts. Previous interpretations have suggested that the bone may have functioned as a lunar or time-reckoning device \parencite{ICOMOS2010, Marshack1971}, or may reflect a base-12 numerical system \parencite{Pletser1999}. Despite these hypotheses, its precise purpose remains debated \parencite{Joseph2011, Rudman2007}.

The starting point of this study is an independent observation of groupings to 30 within Columns~G and~D. Each column contains four numbers totalling 60, and each can be grouped into two pairs summing to 30. This led to the working hypothesis that the numerical arrangement might reflect mathematical intention, which warrants a thorough investigation of the entire set of numbers. In Column~G, the pairs to 30 occupy the outer and inner positions, whereas in Column~D they are crossed, giving the appearance of a difference in positional pairing. Following an exploratory adjustment to ``correct'' this apparent difference in positional pairing, further patterns were observed that span all numerical values in Columns~G and~D, which raised the further hypothesis that the adjustment ``restored'' an underlying arrangement. The present study therefore investigates the patterns and the deviations from them, documents each adjustment together with the patterns that result, and examines whether these findings could provide a basis for demonstrating mathematical intentionality in the numerical arrangement of the Ishango Bone.

The analysis is exploratory: the structural properties and adjustments were developed by examining the published transcription and illustrations of the artefact. Two quantitative analyses, both scoring five specified structural properties, characterise the proposed structure within specified models: a configuration-specific Monte Carlo randomisation, and a complete sensitivity analysis over 64 combinations of readings and adjustments. Their implications are discussed in Section~\ref{sec:stateva} and Section~\ref{sec:discussion}. As a secondary hypothesis, the arithmetic relationships suggest that the bone may have served as a portable reference for demonstrating numerical ideas, for example by arranging stones or other physical markers in a teaching or storytelling setting (Section~\ref{sec:storytelling}).

The sequence of the analysis is central to its interpretation. Before any adjustment, only one pattern spanning all numbers in a column pair was found: the four subgroupings to 30 in Columns G and D. A further relationship (Figure~\ref{fig:Figure11}) already held in the documented arrangement but became apparent only after the exchange in Column~D. Further patterns in Column M, and between Column M and Columns G and D, did not resolve. The first adjustment, the exchange in Column D, was made for the stated reason, namely to align the positional pairing of Column D with that of Column G. Several additional patterns involving all entries of Columns G and D were observed only after this exchange. The analogous exchange in Column M and the unit transfer were then made for stated reasons and were followed by cross-column relationships involving all sixteen entries. The patterns therefore emerged after the adjustments, not before them, and their sequence is documented here as an exploratory finding.

Section~\ref{sec:methods} describes the approach, Sections~\ref{sec:analysisGandD} to~\ref{sec:synthesis} develop the reconstruction, Section~\ref{sec:stateva} presents the quantitative evaluation, and Sections~\ref{sec:storytelling} and~\ref{sec:discussion} discuss the demonstrative interpretation and the limitations.

\vspace{1cm}

\begin{figure}[!t]
    \centering
    \sloppy
    \includegraphics[width=0.95\linewidth]{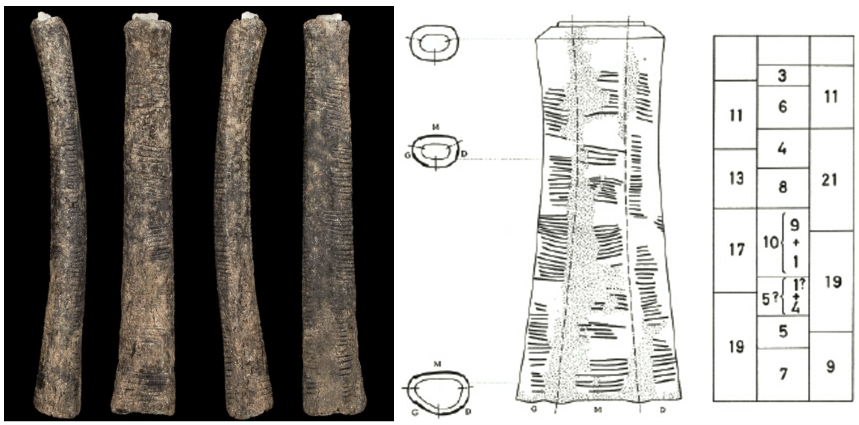}
    \caption{\raggedright Notched bone from Ishango. \\ Left: photograph of the Ishango Bone, © RBINS \parencite{FirstIshangoBone} \\ Right: drawing of the bone by Jean de Heinzelin de Braucourt, © RBINS, from de Heinzelin de Braucourt (1957). Both images reproduced with permission.}
    \label{fig:Figure1}
\end{figure}

\vspace{-0.95cm}

\medskip

This paper is a revised version of an earlier preprint by the author \parencite{baur2025decoding}. This version adds the quantitative analyses (Section~\ref{sec:stateva}) and further patterns.

\section{Methods}
\label{sec:methods}

This study examines de Heinzelin's published transcription and illustrations \parencite[pp.~65--67]{deHeinzelin1957}. The approach is exploratory, using observed relationships among the notch counts to develop hypotheses \parencite{tukey1977}.

The analysis begins with the observation that Columns G and D each total 60 and can be partitioned into two pairs summing to 30. In Column G, these pairs occupy the outer and inner positions. Exchanging the first two entries of Column D (11 and 21) produces the same positional pairing pattern.

Column~M was represented as four consecutive pairs, each aligned with one G--D row; this correspondence is an analytical assumption. The first two M pairs were then exchanged. Under the provisional damaged-pair reading $(9, 4)$, transferring one unit from $5$ to $3$ makes all four M-pair products equal their corresponding G--D row sums, while preserving both the number of entries and their total. Sections~\ref{sec:analysisGandD} and \ref{sec:analysisMGD} present the observations motivating these adjustments.

The uncertain groups Me and Mf in Column M introduce reading uncertainty. De Heinzelin described nine clear marks and an additional damaged mark in Me, and four clear marks and an additional doubtful mark in Mf; his preferred reconstruction was $(10,5)$ \parencite[pp.~66--67]{deHeinzelin1957}. The present analysis also considers the provisional reading $(9,4)$. These two groups are treated as a pair for the proposed row correspondence; this pairing is an analytical assumption.

The reference configurations with documented positions use the readings $(9,4)$ and $(10,5)$ for the damaged pair in Column~M. The sensitivity analysis varies six independently combinable binary choices (each with two alternatives): the two entries of the damaged pair, which reflect reading uncertainty, and four analytical adjustments, namely the exchange in Column~D, the exchange of the first two M pairs, and the unit transfer ($5\rightarrow4$ and $3\rightarrow4$), which are evaluated separately as well as together. 

Both quantitative analyses in this study score an arrangement by how many of five specified structural properties hold at the same time. In the Monte Carlo randomisation, the sixteen numerical entries of each of the three reference configurations were randomly redistributed across the same labelled three-column structure, and the observed score was compared with the scores of 1,000,000 shuffled arrangements per configuration. In the sensitivity analysis, all 64 permitted combinations of readings and adjustments were evaluated once each, without random sampling.

\section{Analysis of Numerical Arrangements in Columns G and D}
\label{sec:analysisGandD}

Jean de Heinzelin de Braucourt documented the grouped notches on the Ishango Bone, which are arranged in three columns. He labelled them M (\textit{milieu}), G (\textit{gauche}) and D (\textit{droite}), and observed that the values in Columns G and D each sum to 60, and that Column~G consists of the prime numbers between 10 and 20. In his schematic illustration, Column M is shown as the middle column \parencite[pp.~65--67]{deHeinzelin1957}. An interactive 3D model of the bone can be explored online via the Royal Belgian Institute of Natural Sciences, where the bone is housed \parencite{FirstIshangoBone}. On the physical object, however, Column M lies on the reverse side of the slightly curved bone, whereas Columns G and D appear on the front. This analysis therefore begins with Columns G and D before turning to Column M.

According to de Heinzelin, both Columns G and D contain four numbers, with 9 being the lowest and 21 being the largest in this set of numbers, as such \parencite[p.~67]{deHeinzelin1957}. Both columns contain only odd numbers. Column~G contains only prime numbers, while Column~D contains two primes, 11 and 19, and two odd composites, 9 and 21, as shown in Figure~\ref{fig:Figure2}.

\begin{figure}[H]
\centering
\begin{tikzpicture}[scale=1]
  \useasboundingbox (-0.5,-1) rectangle (4.5,9.5);
  \draw (0,0) rectangle (4,8);
  \draw (2,0) -- (2,8);
  \draw (0,2) -- (4,2);
  \draw (0,4) -- (4,4);
  \draw (0,6) -- (4,6);
  \node at (1,8.5) {\Large G};
  \node at (3,8.5) {\Large D};
  \node at (1,7) {11};
  \node at (1,5) {13};
  \node at (1,3) {17};
  \node at (1,1) {19};
  \node at (3,7) {11};
  \node at (3,5) {\textcolor{green!60!black}{21}};
  \node at (3,3) {19};
  \node at (3,1) {\textcolor{green!60!black}{9}};
  \node at (1,-0.5) {60};
  \node at (3,-0.5) {60};
\end{tikzpicture}
\caption{All numbers used within Columns G and D. Prime numbers are coloured black, and the two odd composite numbers, 9 and 21, are coloured green.}
\label{fig:Figure2}
\end{figure}

When all integers from 9 to 21 are laid out (odd primes: 11, 13, 17, 19; odd composites: 9, 15, 21; even numbers: 10, 12, 14, 16, 18, 20), as shown in Figure~\ref{fig:Figure3}, it becomes apparent that Columns~G and~D contain every prime number in that range. No even numbers are used, and every odd number in the range appears except 15. The number 15, however, is both the median and the arithmetic mean of the six distinct values present in Columns~G and~D ($9+11+13+17+19+21=90$; $90\div6=15$), as shown in Figure~\ref{fig:Figure4}. If the doubled values 11 and 19 are included, the mean of all eight entries likewise remains 15 ($120\div8=15$).

\begin{figure}[H]
\centering
\begin{tikzpicture}[scale=1]
  \useasboundingbox (-0.5,-1) rectangle (4.5,9.5);
  \draw (0,0) rectangle (4,8);
  \draw (2,0) -- (2,8);
  \draw (0,2) -- (4,2);
  \draw (0,4) -- (4,4);
  \draw (0,6) -- (4,6);
  \node at (1,8.5) {\Large G};
  \node at (3,8.5) {\Large D};
  \node at (1,7) {11};
  \node at (1,5) {13};
  \node at (1,3) {17};
  \node at (1,1) {19};
  \node at (3,7) {11};
  \node at (3,5) {\textcolor{green!60!black}{21}};
  \node at (3,3) {19};
  \node at (3,1) {\textcolor{green!60!black}{9}};
  \node at (1,-0.5) {60};
  \node at (3,-0.5) {60};
\end{tikzpicture}
\vspace{1em}
\begin{minipage}{0.9\linewidth}
\centering
\small
Prime, \textcolor{green!60!black}{Odd composite}, and
\textcolor{pink!80!red}{Even} numbers from 9 to 21:\\[0.3em]
\textcolor{green!60!black}{9}, \textcolor{pink!80!red}{10}, 11, \textcolor{pink!80!red}{12}, 13, \textcolor{pink!80!red}{14}, \textcolor{green!60!black}{15}, \textcolor{pink!80!red}{16}, 17, \textcolor{pink!80!red}{18}, 19, \textcolor{pink!80!red}{20}, \textcolor{green!60!black}{21}
\end{minipage}
\caption{All numbers ranging from 9 to 21 are colour-coded by mathematical property. All primes in Columns G and D are included, as well as all odd numbers except 15. No even numbers are used.}
\label{fig:Figure3}
\end{figure}

When all odd numbers from 9 to 21 are arranged in sequence, as shown in Figure~\ref{fig:Figure4}, they form an arithmetic progression with a common difference of 2 ($9,11,13,15,17,19,21$). The two outer numbers sum to 30 ($9+21$), and so do the middle pair ($11+19$) and the inner pair ($13+17$). Columns~G and~D contain all numbers and pairs of this sequence except 15, with the middle pair, 11 and 19, occurring in both columns.

\vspace{1cm}

\begin{figure}[H]
\centering
\begin{tikzpicture}[scale=1]
  \node at (0,0) {\textcolor{green!60!black}{9}};
  \node at (1.2,0) {11};
  \node at (2.4,0) {13};
  \node at (3.6,0) {\textcolor{green!60!black}{15}};
  \node at (4.8,0) {17};
  \node at (6.0,0) {19};
  \node at (7.2,0) {\textcolor{green!60!black}{21}};
  \draw (2.4,-0.3) -- (2.4,-0.7) -- (4.8,-0.7) -- (4.8,-0.3);
  \node at (3.6,-1.0) {\small 30};
  \draw (1.2,-0.3) -- (1.2,-1.4) -- (6.0,-1.4) -- (6.0,-0.3);
  \node at (3.6,-1.7) {\small 30};
  \draw (0,-0.3) -- (0,-2.1) -- (7.2,-2.1) -- (7.2,-0.3);
  \node at (3.6,-2.4) {\small 30};
\end{tikzpicture}
\vspace{1em}
\caption{All odd numbers ranging from 9 to 21 are lined up. The prime numbers are coloured black, and the odd composite numbers are coloured green. Three distinct groupings can be formed around the centre number 15: the inner ($13+17$), middle ($11+19$), and outer grouping ($9+21$), each summing to 30.}
\label{fig:Figure4}
\end{figure}

When this grouping rule is applied to the numbers in Columns~G and~D, as shown in Figure~\ref{fig:Figure5}, each column can be divided into two pairs summing to 30. In Column~G, these are the outer and inner pairs ($11+19$ and $13+17$), whereas in Column~D they are crossed ($11+19$ and $21+9$), so that the pairs share the sum of 30 but not the same positional arrangement.

\begin{figure}[H]
\centering
\begin{tikzpicture}[scale=1]
  \draw (0,0) rectangle (4,8);
  \draw (2,0) -- (2,8);
  \draw (0,2) -- (4,2);
  \draw (0,4) -- (4,4);
  \draw (0,6) -- (4,6);
  \node at (1,8.5) {\Large G};
  \node at (3,8.5) {\Large D};
  \node at (1,7) {11};
  \node at (1,5) {13};
  \node at (1,3) {17};
  \node at (1,1) {19};
  \node at (3,7) {11};
  \node at (3,5) {\textcolor{green!60!black}{21}};
  \node at (3,3) {19};
  \node at (3,1) {\textcolor{green!60!black}{9}};
  \node at (1,-0.5) {60};
  \node at (3,-0.5) {60};
  \draw (0,7) -- (-1.0,7) -- (-1.0,1) -- (0,1);
  \node[rotate=90] at (-1.3,4) {30};
  \draw (0,5) -- (-0.5,5) -- (-0.5,3) -- (0,3);
  \node[rotate=90] at (-0.75,4) {30};
  \draw (4,5) -- (4.9,5) -- (4.9,1) -- (4,1);
  \node[rotate=270] at (5.15,3) {30};
\draw (4,7) -- (4.4,7) -- (4.4,3) -- (4,3);
\node[rotate=270] at (4.65,6) {30};
\end{tikzpicture}
\caption{Four distinct groupings can be identified within Columns G and D, with each column containing two groupings that each sum to 30, following the same grouping rule as shown in Figure~\ref{fig:Figure4}. The grouping, $11 + 19$, appears twice, once in each column. The groupings in Column D are crossed positional pairings, while they are structured in Column G.}
\label{fig:Figure5}
\end{figure}

In summary, Columns~G and~D together contain every odd number between 9 and 21 except 15. Each column sums to 60 and can be divided into two pairs that sum to 30. In Column~G these pairs occupy the outer and inner positions, whereas in Column~D they are crossed, a deviation that is examined in the following section.

\subsection*{Aligning the Symmetry in Column D}

Columns~G and~D contain only odd numbers, each totals 60, and each can be partitioned into two pairs summing to 30. The present study therefore takes intentional numerical organisation as a working hypothesis. Relative to the outer--inner arrangement in Column~G, the complementary pairs in Column~D appear positionally misaligned, and a transmission error is one possible explanation.

Models of cultural transmission show how copying errors can generate variation in material culture \parencite{eerkens2005}. A transmission explanation prompted the initial exchange in Column~D.

The numbers 11 and 21 in Column~D were then exchanged (Figure~\ref{fig:Figure6}).

\begin{figure}[H]
\centering
\begin{tikzpicture}[scale=1]
\useasboundingbox (-0.5,-1) rectangle (5.0,9.5);
\draw (0,0) rectangle (4,8);
\draw (2,0) -- (2,8);
\draw (0,2) -- (4,2);
\draw (0,4) -- (4,4);
\draw (0,6) -- (4,6);
\node at (1,8.5) {\Large G};
\node at (3,8.5) {\Large D};
\node at (1,7) {11};
\node at (1,5) {13};
\node at (1,3) {17};
\node at (1,1) {19};
\node at (3,7) {11};
\node at (3,5) {21};
\node at (3,3) {19};
\node at (3,1) {9};
\node at (1,-0.6) {60};
\node at (3,-0.6) {60};
\draw[green!60!black]
  (4.15,6.9) -- (4.5,6.9)
  -- (4.5,4.9)
  -- (4.15,4.9);
\fill[green!60!black]
  (4.02,6.9) -- (4.15,7.0) -- (4.15,6.8);
\fill[green!60!black]
  (4.02,4.9) -- (4.15,5.0) -- (4.15,4.8);
\end{tikzpicture}
\caption{The first two entries of Column D, 11 and 21, are exchanged to align its complementary pairings with the outer--inner arrangement observed in Column G. The additional relationships discussed below were noticed after this exchange.}
\label{fig:Figure6}
\end{figure}

After this exchange, Column~D follows the same outer--inner grouping as Column~G: the outer pair ($21+9$) and the inner pair ($11+19$) each sum to 30, as shown in Figure~\ref{fig:Figure7}.

\begin{figure}[H]
\centering
\begin{tikzpicture}[scale=1]
  \draw (0,0) rectangle (4,8);
  \draw (2,0) -- (2,8);
  \draw (0,2) -- (4,2);
  \draw (0,4) -- (4,4);
  \draw (0,6) -- (4,6);
  \node at (1,8.5) {\Large G};
  \node at (3,8.5) {\Large D};
  \node at (1,7) {11};
  \node at (1,5) {13};
  \node at (1,3) {17};
  \node at (1,1) {19};
  \node at (3,7) {21};
  \node at (3,5) {11};
  \node at (3,3) {19};
  \node at (3,1) {9};
  \node at (1,-0.5) {60};
  \node at (3,-0.5) {60};
  \draw (0,7) -- (-1.0,7) -- (-1.0,1) -- (0,1);
  \node[rotate=90] at (-1.3,4) {30};
  \draw (0,5) -- (-0.5,5) -- (-0.5,3) -- (0,3);
  \node[rotate=90] at (-0.75,4) {30};
  \draw (4,7) -- (4.9,7) -- (4.9,1) -- (4,1);
  \node[rotate=270] at (5.15,4) {30};
  \draw (4,5) -- (4.4,5) -- (4.4,3) -- (4,3);
  \node[rotate=270] at (4.65,4) {30};
\end{tikzpicture}
\caption{Upon exchanging the position of numbers 11 and 21 in Column D, the resulting arrangement appears to align with the symmetrical groupings seen in Column G, and the grouping in Column D no longer crossed.}
\label{fig:Figure7}
\end{figure}

After this exchange, the adjusted arrangement was examined further, and several patterns across Columns~G and~D could be observed.

One pattern, illustrated in Figure~\ref{fig:Figure8}, concerns the two numbers that occur in both columns, 11 and 19. Each occupies the two diagonal cells of one block: 11 appears only in the upper block and 19 only in the lower block, with one copy in each column and each row of its block.

\begin{figure}[H]
\centering
\begin{tikzpicture}[scale=1]
\useasboundingbox (-0.5,-1) rectangle (5.0,9.5);
\fill[warmyellow] (0,6) rectangle (2,8); 
\fill[warmyellow] (2,4) rectangle (4,6); 
\fill[sagegreen]  (2,2) rectangle (4,4); 
\fill[sagegreen]  (0,0) rectangle (2,2); 
\draw (0,0) rectangle (4,8);
\draw (2,0) -- (2,8);
\draw (0,2) -- (4,2);
\draw (0,4) -- (4,4);
\draw (0,6) -- (4,6);
\node at (1,8.5) {\Large G};
\node at (3,8.5) {\Large D};
\node at (1,7) {11};
\node at (1,5) {13};
\node at (1,3) {17};
\node at (1,1) {19};
\node at (3,7) {21};
\node at (3,5) {11};
\node at (3,3) {19};
\node at (3,1) {9};
\end{tikzpicture}
\caption{The rearrangement of numbers 11 and 21 in Column D resulted in an additional pattern: the doubled prime numbers 11 and 19 are now placed diagonally to each other across Columns G and D. Number 11 appears twice in the upper quadrant and in each row and column once. Number 19 appears twice in the lower quadrant and in each row and column once.}
\label{fig:Figure8}
\end{figure}

A further pattern is shown in Figure~\ref{fig:Figure9}. Each of the four unique numbers can be added to the doubled number in its row and to the doubled number in its column, and the two sums are identical for every unique number: $13+11=24$, $21+11=32$, $17+19=36$ and $9+19=28$. Together, these give the sums 24, 28, 32 and 36.

\begin{figure}[H]
\centering
\begin{tikzpicture}[scale=1]
\useasboundingbox (-0.5,-1) rectangle (8.0,9.5);
\fill[warmyellow] (0,6) rectangle (2,8); 
\fill[warmyellow] (2,4) rectangle (4,6); 
\fill[sagegreen]  (2,2) rectangle (4,4); 
\fill[sagegreen]  (0,0) rectangle (2,2); 
\draw (0,0) rectangle (4,8);
\draw (2,0) -- (2,8);
\draw (0,2) -- (4,2);
\draw (0,4) -- (4,4);
\draw (0,6) -- (4,6);
\node at (1,8.5) {\Large G};
\node at (3,8.5) {\Large D};
\node at (1,7) {11};
\node at (1,5) {13};
\node at (1,3) {17};
\node at (1,1) {19};
\node at (3,7) {21};
\node at (3,5) {11};
\node at (3,3) {19};
\node at (3,1) {9};
\draw[->, thick] (2.6,7) -- (1.4,7);
\draw[->, thick] (3,6.6) -- (3,5.4);
\draw[->, thick] (1.4,5) -- (2.6,5);
\draw[->, thick] (1,5.4) -- (1,6.6);
\draw[->, thick] (1.4,3) -- (2.6,3);
\draw[->, thick] (1,2.6) -- (1,1.4);
\draw[->, thick] (2.6,1) -- (1.4,1);
\draw[->, thick] (3,1.4) -- (3,2.6);
\node at (5.5,7) {$11+21=32$};
\node at (5.5,5) {$13+11=24$};
\node at (5.5,3) {$17+19=36$};
\node at (5.5,1) {$19+9=28$};
\end{tikzpicture}
\caption{Each unique number in Columns G and D can be added to the adjacent doubled number in its row and to the doubled number in its column, and both sums are identical. For $21$ (Column D), both sums are $21 + 11 = 32$. For $13$ (Column G), both are $13 + 11 = 24$. For $17$ (Column G), both are $17 + 19 = 36$. For $9$ (Column D), both are $9 + 19 = 28$.}
\label{fig:Figure9}
\end{figure}

It was also observed, as shown in Figure~\ref{fig:Figure10}, that the outer groupings of both columns ($11+19$ in Column~G and $21+9$ in Column~D) together sum to 60, which equals the total of the inner groupings ($13+17$ in Column~G and $11+19$ in Column~D). In terms of rows, the outer rows (1 and 4) and the inner rows (2 and 3) each total 60.

\begin{figure}[H]
\centering
\begin{tikzpicture}[scale=1]
\useasboundingbox (-0.5,-1) rectangle (8.5,9.5);
  \draw (0,0) rectangle (4,8);
  \draw (2,0) -- (2,8);
  \draw (0,2) -- (4,2);
  \draw (0,4) -- (4,4);
  \draw (0,6) -- (4,6);
  \node at (1,8.5) {\Large G};
  \node at (3,8.5) {\Large D};
  \node at (1,7) {11};
  \node at (1,5) {13};
  \node at (1,3) {17};
  \node at (1,1) {19};
  \node at (3,7) {21};
  \node at (3,5) {11};
  \node at (3,3) {19};
  \node at (3,1) {9};
  \node at (1,-0.5) {60};
  \node at (3,-0.5) {60};
  \draw (0,7) -- (-1.0,7) -- (-1.0,1) -- (0,1);
  \node[rotate=90] at (-1.3,4) {30};
  \draw (0,5) -- (-0.5,5) -- (-0.5,3) -- (0,3);
  \node[rotate=90] at (-0.75,4) {30};
  \draw (4,7) -- (4.9,7) -- (4.9,1) -- (4,1);
  \node[rotate=270] at (5.15,4) {30};
  \draw (4,5) -- (4.4,5) -- (4.4,3) -- (4,3);
  \node[rotate=270] at (4.65,4) {30};
  \node at (6.5,7) {$11+21=32$};
  \node at (6.5,5) {$13+11=24$};
  \node at (6.5,3) {$17+19=36$};
  \node at (6.5,1) {$19+9=28$};
  \draw (8.1,7) -- (9.0,7) -- (9.0,1) -- (8.1,1);
  \node[rotate=270] at (9.25,4) {60};
  \draw (8.1,5) -- (8.5,5) -- (8.5,3) -- (8.1,3);
  \node[rotate=270] at (8.75,4) {60};
\end{tikzpicture}
\caption{The position exchange of numbers 11 and 21 in Column D also resulted in an extended grouping pattern where the upper and lower rows, as well as the two middle rows, combined, summing to 60 in each case.}
\label{fig:Figure10}
\end{figure}

A further relationship concerns the upper block (rows 1 and 2) and the lower block (rows 3 and 4). Their totals are 56 and 64, respectively; subtracting one occurrence of the duplicated value in each block gives $56-11=45$ and $64-19=45$ (Figure~\ref{fig:Figure11}). However, this relationship already holds in the documented arrangement, since exchanging the first two entries of Column~D only changes positions within the upper block and leaves both block totals
unchanged. It is therefore not an additional equality produced by the exchange.

\begin{figure}[H]
\centering
\begin{tikzpicture}[scale=1]
\useasboundingbox (-0.5,-1) rectangle (7.0,9.5);
  \fill[warmyellow] (0,6) rectangle (2,8); 
  \fill[warmyellow] (2,4) rectangle (4,6); 
  \fill[sagegreen]  (2,2) rectangle (4,4); 
  \fill[sagegreen]  (0,0) rectangle (2,2); 
  \draw (0,0) rectangle (4,8);
  \draw (2,0) -- (2,8);
  \draw (0,2) -- (4,2);
  \draw (0,4) -- (4,4);
  \draw (0,6) -- (4,6);
  \node at (1,8.5) {\Large G};
  \node at (3,8.5) {\Large D};
  \node at (1,7) {11};
  \node at (1,5) {13};
  \node at (1,3) {17};
  \node at (1,1) {19};
  \node at (3,7) {21};
  \node at (3,5) {11};
  \node at (3,3) {19};
  \node at (3,1) {9};
  \draw[decorate, decoration={brace, amplitude=12pt}, line width=0.8pt] 
    (4,7) -- (4,5);
  \node[right] at (4.8,6) {$56 - 11 = 45$};
  \draw[decorate, decoration={brace, amplitude=12pt}, line width=0.8pt] 
    (4,3) -- (4,1);
  \node[right] at (4.8,2) {$64 - 19 = 45$};
\end{tikzpicture}
\caption{When the upper and lower quadrants are summed independently and one of their duplicated numbers is subtracted, each block yields the same total of 45.}
\label{fig:Figure11}
\end{figure}

In summary, this single positional exchange was followed by further patterns involving all values in Columns~G and~D. Together, these patterns are consistent with a mathematical organisation of the numbers. Columns~G and~D, however, account for only two of the three columns, and a fuller picture requires examining Column~M, which is addressed in the next section.

\section{Structural Analysis of Column M and Its Relationship to Columns G and D}
\label{sec:analysisMGD}

The eight numerical values of Column~M, grouped into four pairs, were initially examined in relation to the four rows of Columns~G and~D, with each pair tentatively aligned to one row (Figure~\ref{fig:Figure12}).

\begin{figure}[H]
\centering
\begin{tikzpicture}[scale=1]
\useasboundingbox (-0.5,-1) rectangle (8.0,9.5);
  \draw (0,0) rectangle (4,8);
  \draw (2,0) -- (2,8);
  \draw (0,2) -- (4,2);
  \draw (0,4) -- (4,4);
  \draw (0,6) -- (4,6);
  \node at (1,8.5) {\Large G};
  \node at (3,8.5) {\Large D};
  \node at (1,7) {11};
  \node at (1,5) {13};
  \node at (1,3) {17};
  \node at (1,1) {19};
  \node at (3,7) {21};
  \node at (3,5) {11};
  \node at (3,3) {19};
  \node at (3,1) {9};
  \draw (5.5,0) rectangle (7.5,8);
  \draw (5.5,2) -- (7.5,2);
  \draw (5.5,4) -- (7.5,4);
  \draw (5.5,6) -- (7.5,6);
  \node at (6.5,8.5) {\Large M};
  \node at (6.5,7.2) {3};
  \node at (6.5,6.8) {6};
  \node at (6.5,5.2) {4};
  \node at (6.5,4.8) {8};
  \node at (6.5,3.2) {9(+1?)};
  \node at (6.5,2.8) {4(+1?)};
  \node at (6.5,1.2) {5};
  \node at (6.5,0.8) {7};
\end{tikzpicture}
\caption{Including Column M in the further analysis, the eight numbers in Column M are arranged in pairs, each corresponding by row to the four numbers in Columns G and D.}
\label{fig:Figure12}
\end{figure}

The numerical values in Column~M were first examined without adjustment, but no structural pattern applying to all values was found, either between Column~M and Columns~G and~D or within Column~M itself. The exchange in Column~D had produced consolidated patterns across Columns~G
and~D. To test whether the same would occur in Column~M, the first two pairs in Column~M were exchanged analogously, following the same positional logic. This step is an instance of abductive reasoning, in which a possible explanation is inferred from observed relationships
\parencite{peirce1903}. Because the exchange was developed through sequential examination of the same artefact, it is not equivalent to an independently tested prediction \parencite{popper1959}.

\begin{figure}[H]
\centering
\begin{tikzpicture}[scale=1]
\useasboundingbox (-0.5,-1) rectangle (9.0,9.5);
  \draw (0,0) rectangle (4,8);
  \draw (2,0) -- (2,8);
  \draw (0,2) -- (4,2);
  \draw (0,4) -- (4,4);
  \draw (0,6) -- (4,6);
  \node at (1,8.5) {\Large G};
  \node at (3,8.5) {\Large D};
  \node at (1,7) {11};
  \node at (1,5) {13};
  \node at (1,3) {17};
  \node at (1,1) {19};
  \node at (3,7) {21};
  \node at (3,5) {11};
  \node at (3,3) {19};
  \node at (3,1) {9};
  \draw (5.5,0) rectangle (7.5,8);
  \draw (5.5,2) -- (7.5,2);
  \draw (5.5,4) -- (7.5,4);
  \draw (5.5,6) -- (7.5,6);
  \node at (6.5,8.5) {\Large M};
  \node at (6.5,7.2) {3};
  \node at (6.5,6.8) {6};
  \node at (6.5,5.2) {4};
  \node at (6.5,4.8) {8};
  \node at (6.5,3.2) {9(+1?)};
  \node at (6.5,2.8) {4(+1?)};
  \node at (6.5,1.2) {5};
  \node at (6.5,0.8) {7};
\draw[green!60!black]
    (7.65,7.0) -- (7.95,7.0)
    -- (7.95,5.0)
    -- (7.65,5.0);
\fill[green!60!black]
    (7.5,7.0) -- (7.65,7.08) -- (7.65,6.92);
\fill[green!60!black]
    (7.5,5.0) -- (7.65,5.08) -- (7.65,4.92);
\draw[black]
    (4.15,7.0) -- (4.45,7.0)
    -- (4.45,5.0)
    -- (4.15,5.0);
\fill[black]
    (4.0,7.0) -- (4.15,7.08) -- (4.15,6.92);
\fill[black]
    (4.0,5.0) -- (4.15,5.08) -- (4.15,4.92);
\end{tikzpicture}
\caption{The number pairs (3, 6) and (4, 8) in Column M were exchanged, mirroring the position exchange applied to Column D, to test whether the positional deviation in Column D reflected a structural feature of the design rather than a transmission error, and whether additional mathematical relationships would emerge as they did following the Column D exchange.}
\label{fig:Figure13}
\end{figure}

After the first two pairs in Column M were exchanged, a possible multiplicative relationship became apparent, as shown in Figure~\ref{fig:Figure14}. Each row of Columns G and D has a sum, and each row of Column M has a pair of numbers that can be multiplied. In the first and third rows, the product of the pair equals the row sum exactly: the pairs $(4,8)$ and $(9,4)$ give $32$ and $36$.

In each of the other two rows, changing one factor by one unit would make its product equal the corresponding G--D row sum. In the second row, the pair $(3,6)$ gives $18$, but the row sum is $24$; the $3$ would need to be $4$. In the fourth row, the pair $(5,7)$ gives $35$, but the row sum is $28$; the $5$ would need to be $4$. Moving one unit from the $5$ to the $3$ fixes both rows at once. The result is that every row sum in Columns $G$ and $D$ equals
the product of the pair in Column M:

\begin{equation*}
  11 + 21 = 4 \times 8, \qquad 13 + 11 = 4 \times 6, \qquad
  17 + 19 = 9 \times 4, \qquad 19 + 9 = 4 \times 7.
\end{equation*}

The required factor changes are equal and opposite: \(3\rightarrow4\) and \(5\rightarrow4\). A single unit transfer therefore produces both equalities while preserving the combined total.

\begin{figure}[H]
\centering
\begin{tikzpicture}[scale=1]
\useasboundingbox (-0.5,-1) rectangle (12.0,9.5);
  \draw (0,0) rectangle (4,8);
  \draw (2,0) -- (2,8);
  \draw (0,2) -- (4,2);
  \draw (0,4) -- (4,4);
  \draw (0,6) -- (4,6);
  \node at (1,8.5) {\Large G};
  \node at (3,8.5) {\Large D};
  \node at (1,7) {11};
  \node at (1,5) {13};
  \node at (1,3) {17};
  \node at (1,1) {19};
  \node at (3,7) {21};
  \node at (3,5) {11};
  \node at (3,3) {19};
  \node at (3,1) {9};
  \node at (5.5,7) {$11+21=$\textcolor{green!60!black}{$32$}};
  \node at (5.5,5) {$13+11=$\textcolor{green!60!black}{$24$}};
  \node at (5.5,3) {$17+19=$\textcolor{green!60!black}{$36$}};
  \node at (5.5,1) {$19+9=$\textcolor{green!60!black}{$28$}};
  \draw (7.0,0) rectangle (9.0,8);
  \draw (7.0,2) -- (9.0,2);
  \draw (7.0,4) -- (9.0,4);
  \draw (7.0,6) -- (9.0,6);
  \node at (8.0,8.5) {\Large M};
  \node at (8.0,7.2) {4};
  \node at (8.0,6.8) {8};
  \node at (8.0,5.2) {\textcolor{orange}{3}};
  \node at (8.0,4.8) {6};
  \node at (8.0,3.2) {9(+1?)};
  \node at (8.0,2.8) {4(+1?)};
  \node at (8.0,1.2) {\textcolor{orange}{5}};
  \node at (8.0,0.8) {7};
  \node at (10.5,7) {$4\times8=$\textcolor{green!60!black}{$32$}};
  \node at (10.5,5) {$4\times6=$\textcolor{green!60!black}{$24$}};
  \node at (10.5,3) {$9\times4=$\textcolor{green!60!black}{$36$}};
  \node at (10.5,1) {$4\times7=$\textcolor{green!60!black}{$28$}};
\end{tikzpicture}
\caption{An apparent correspondence between the G--D row sums and the products of the aligned M pairs emerges following the positional exchanges. Two pairs match exactly if the entries at positions 9 or 10 and 4 or 5 are interpreted as 9 and 4, respectively. For the remaining two pairs, interpreting the entries 3 and 5 as 4 would produce the equalities shown on the right.}
\label{fig:Figure14}
\end{figure}

A single-unit transfer from $5$ to $3$ therefore resolves both discrepancies, changing the entries to $4$ and $4$ while preserving their combined total (Figure~\ref{fig:Figure15}). Under these fixed assumptions, no other single-unit transfer resolves both rows.

\begin{figure}[H]
\centering
\begin{tikzpicture}[scale=1]
\useasboundingbox (-1.0,-1.5) rectangle (10.5,9.5);
\draw (0,0) rectangle (2,8);
\draw (0,2) -- (2,2);
\draw (0,4) -- (2,4);
\draw (0,6) -- (2,6);
\node at (1,8.5) {\Large M};
\node at (1,7.2) {4};
\node at (1,6.8) {8};
\node at (1,5.2) {\textcolor{orange}{3(+1)}};
\node at (1,4.8) {6};
\node at (1,3.2) {9(+1?)};
\node at (1,2.8) {4(+1?)};
\node at (1,1.2) {\textcolor{orange}{5(-1)}};
\node at (1,0.8) {7};
\node at (1,-0.8) {Before unit transfer};
\draw (5,0) rectangle (7,8);
\draw (5,2) -- (7,2);
\draw (5,4) -- (7,4);
\draw (5,6) -- (7,6);
\node at (6,8.5) {\Large M};
\node at (6,7.2) {4};
\node at (6,6.8) {8};
\node at (6,5.2) {\textcolor{green!60!black}{4}};
\node at (6,4.8) {6};
\node at (6,3.2) {9};
\node at (6,2.8) {4};
\node at (6,1.2) {\textcolor{green!60!black}{4}};
\node at (6,0.8) {7};
\node at (6,-0.8) {Adjusted};
\draw[green!60!black, dashed]
  (2.15,5.2) -- (2.45,5.2)
  -- (2.45,1.2)
  -- (2.15,1.2);
\fill[green!60!black]
  (2.0,5.2) -- (2.15,5.28) -- (2.15,5.12);
\fill[green!60!black]
  (2.0,1.2) -- (2.15,1.28) -- (2.15,1.12);
\node[rotate=270] at (2.75,3.2) {\small\textcolor{green!60!black}{$\pm$1}};
\end{tikzpicture}
\caption{For further evaluation, one unit was reallocated from the value 5 to the value 3 in Column M. This single unit reallocation resolves both deviating values to 4.}
\label{fig:Figure15}
\end{figure}

\subsection*{Identification of a New Pattern and Mathematical Features within Column M}

Further observation of the adjusted Column~M values showed that, after setting aside one factor of 4 from each M pair, the remaining values, referred to here as the primary values, are $(8,6,9,7)$ in displayed order. When sorted, these form the consecutive sequence $(6,7,8,9)$. Applying the outer--inner grouping used for Columns~G and~D to the primary values in displayed order gives two equal sums: $8+7=15$ and $6+9=15$ (Figure~\ref{fig:Figure16}).

\begin{figure}[H]
\centering
\begin{tikzpicture}[scale=1]
\useasboundingbox (-0.5,-1) rectangle (4.5,9);
\draw (0,0) rectangle (2,8);
\draw (0,2) -- (2,2);
\draw (0,4) -- (2,4);
\draw (0,6) -- (2,6);
\node at (1,8.6) {\Large M};
\node at (1,7) {8};
\node at (1,5) {6};
\node at (1,3) {9};
\node at (1,1) {7};
\draw (2,7) -- (3.0,7)
      -- (3.0,1)
      -- (2,1);
\node[rotate=270] at (3.3,4) {15};
\draw (2,5) -- (2.5,5)
      -- (2.5,3)
      -- (2,3);
\node[rotate=270] at (2.75,4) {15};
\end{tikzpicture}
\caption{A grouping pattern was identified in Column M within the primary numbers, corresponding to the same grouping patterns observed in Columns G and D, where the upper and lower numbers formed the same sum as the two middle numbers in each column. In Column M, the upper and lower number pairs also form two groupings, with each pair summing to 15, collectively yielding a total of 30.}
\label{fig:Figure16}
\end{figure}

The adjusted M pairs thus match the G--D row sums, while their primary values share the outer--inner grouping rule: $8 + 7 = 6 + 9 = 15$. This common sum is also the odd number absent from the interval represented in Columns G and D.

\section{Cross-Column Arithmetic Relationships}
\label{sec:arithmeticalRel}

In each adjusted M pair, the value 4 is termed secondary and the remaining value primary. Figure~\ref{fig:Figure17} applies the secondary value to the primary values by subtraction and addition.

\begin{figure}[H]
\centering
\begin{tikzpicture}[scale=1]
\useasboundingbox (-4,-1) rectangle (6,9.5);
\draw (0,0) rectangle (2,8);
\draw (0,2) -- (2,2);
\draw (0,4) -- (2,4);
\draw (0,6) -- (2,6);
\node at (1,8.6) {\Large M};
\node at (0.7,7) {4};
\node[green!60!black] at (1.3,7) {8};
\node at (0.7,5) {4};
\node[green!60!black] at (1.3,5) {6};
\node at (0.7,3) {4};
\node[green!60!black] at (1.3,3) {9};
\node at (0.7,1) {4};
\node[green!60!black] at (1.3,1) {7};
\node at (-2,8.6) {Subtraction};
\node at (-2,7) {$8-4=$\textcolor{green!60!black}{4}};
\node at (-2,5) {$6-4=$\textcolor{green!60!black}{2}};
\node at (-2,3) {$9-4=$\textcolor{green!60!black}{5}};
\node at (-2,1) {$7-4=$\textcolor{green!60!black}{3}};
\node at (4,8.6) {Addition};
\node at (4,7) {$8+4=$\textcolor{green!60!black}{12}};
\node at (4,5) {$6+4=$\textcolor{green!60!black}{10}};
\node at (4,3) {$9+4=$\textcolor{green!60!black}{13}};
\node at (4,1) {$7+4=$\textcolor{green!60!black}{11}};
\end{tikzpicture}
\caption{The four primary numbers in Column M (6, 7, 8, 9) relate to the supplementary number 4 through both subtraction and addition. Subtracting 4 from each primary number produces the sequence 2, 3, 4, and 5, while adding 4 extends the sequence to 10, 11, 12, and 13.}
\label{fig:Figure17}
\end{figure}

Subtracting $4$ from the primary values $(6, 7, 8, 9)$ produces $(2, 3, 4, 5)$, while adding $4$ produces $(10, 11, 12, 13)$. Together, the three blocks span the integers from $2$ to $13$ without interruption or overlap.

The multiplicative relationship is as follows: in the adjusted arrangement, the product of each M pair equals the sum of the corresponding G--D row (Figure~\ref{fig:Figure18}).

\begin{figure}[H]
\centering
\begin{tikzpicture}[scale=1]
\useasboundingbox (-0.5,-1) rectangle (13,9.5);
\node at (6,9) {\large Multiplication};
\draw (0,0) rectangle (4,8);
\draw (2,0) -- (2,8);
\draw (0,2) -- (4,2);
\draw (0,4) -- (4,4);
\draw (0,6) -- (4,6);
\node at (1,8.5) {\Large G};
\node at (3,8.5) {\Large D};
\node at (1,7) {11};
\node at (1,5) {13};
\node at (1,3) {17};
\node at (1,1) {19};
\node at (3,7) {21};
\node at (3,5) {11};
\node at (3,3) {19};
\node at (3,1) {9};
\node at (6,7) {$11+21=$\textcolor{green!60!black}{32}};
\node at (6,5) {$13+11=$\textcolor{green!60!black}{24}};
\node at (6,3) {$17+19=$\textcolor{green!60!black}{36}};
\node at (6,1) {$19+9=$\textcolor{green!60!black}{28}};
\draw (8,0) rectangle (10,8);
\draw (8,2) -- (10,2);
\draw (8,4) -- (10,4);
\draw (8,6) -- (10,6);
\node at (9,8.5) {\Large M};
\node at (9,7.3) {4};
\node at (9,6.7) {8};
\node at (9,5.3) {4};
\node at (9,4.7) {6};
\node at (9,3.3) {9};
\node at (9,2.7) {4};
\node at (9,1.3) {4};
\node at (9,0.7) {7};
\node at (11.5,7) {\textcolor{green!60!black}{32}$=4\times8$};
\node at (11.5,5) {\textcolor{green!60!black}{24}$=4\times6$};
\node at (11.5,3) {\textcolor{green!60!black}{36}$=4\times9$};
\node at (11.5,1) {\textcolor{green!60!black}{28}$=4\times7$};
\end{tikzpicture}
\caption{In the adjusted configuration, the sum of the G and D entries in each row equals the product of the two entries in the corresponding M pair: $32=4\times8$, $24=4\times6$, $36=4\times9$, and $28=4\times7$.}
\label{fig:Figure18}
\end{figure}

Division gives equivalent relationships that follow directly from the products (Figure~\ref{fig:Figure19}). Each row sum divided by the primary value of its aligned M pair yields the secondary value 4 ($32\div8=4$, $24\div6=4$, $36\div9=4$, $28\div7=4$). Equivalently, dividing each row sum by the secondary value, number 4, yields the primary value of the aligned pair ($32\div4=8$, $24\div4=6$, $36\div4=9$, $28\div4=7$).

\begin{figure}[H]
\centering
\begin{tikzpicture}[scale=1]
\useasboundingbox (-0.5,-1) rectangle (13,9.5);
\node at (6,9) {\large Division};
\draw (0,0) rectangle (4,8);
\draw (2,0) -- (2,8);
\draw (0,2) -- (4,2);
\draw (0,4) -- (4,4);
\draw (0,6) -- (4,6);
\node at (1,8.5) {\Large G};
\node at (3,8.5) {\Large D};
\node at (1,7) {11};
\node at (1,5) {13};
\node at (1,3) {17};
\node at (1,1) {19};
\node at (3,7) {21};
\node at (3,5) {11};
\node at (3,3) {19};
\node at (3,1) {9};
\node at (6,7) {$11+21=$\textcolor{green!60!black}{32}};
\node at (6,5) {$13+11=$\textcolor{green!60!black}{24}};
\node at (6,3) {$17+19=$\textcolor{green!60!black}{36}};
\node at (6,1) {$19+9=$\textcolor{green!60!black}{28}};
\draw (8,0) rectangle (10,8);
\draw (8,2) -- (10,2);
\draw (8,4) -- (10,4);
\draw (8,6) -- (10,6);
\node at (9,8.5) {\Large M};
\node at (9,7.3) {4};
\node at (9,6.7) {8};
\node at (9,5.3) {4};
\node at (9,4.7) {6};
\node at (9,3.3) {9};
\node at (9,2.7) {4};
\node at (9,1.3) {4};
\node at (9,0.7) {7};
\node at (11.5,7) {\textcolor{green!60!black}{32}$\div4=8$};
\node at (11.5,5) {\textcolor{green!60!black}{24}$\div4=6$};
\node at (11.5,3) {\textcolor{green!60!black}{36}$\div4=9$};
\node at (11.5,1) {\textcolor{green!60!black}{28}$\div4=7$};
\end{tikzpicture}
\caption{For each row, the sum of Columns G and D divided by the secondary value, number 4, of its corresponding Column M pair yields the primary value; equivalently, dividing by the primary value yields the secondary.}
\label{fig:Figure19}
\end{figure}

The adjusted configuration contains relationships that can be expressed with all four basic arithmetic operations. Such relationships could be used to demonstrate how different operations connect the same quantities, which provides one motivation for the demonstrative-use hypothesis considered in Section~\ref{sec:storytelling}.

\section{Synthesis: Cross-Column Grouping and Progression Relationships}
\label{sec:synthesis}

In the adjusted configuration, the G--D row sums are $(32,24,36,28)$ in displayed order. When sorted, they form the arithmetic progression $(24,28,32,36)$ with common difference 4. Adding the corresponding M primary values $(8,6,9,7)$ gives $(40,30,45,35)$, which form the arithmetic progression $(30,35,40,45)$ when sorted, with common difference 5, as shown in Figure~\ref{fig:Figure20}.

\begin{figure}[H]
\centering
\begin{tikzpicture}[scale=1]
\useasboundingbox (-0.5,-4) rectangle (13,10.2);
\node at (6,9.8) {\large Progression of 4 and 5};
\draw (0,0) rectangle (4,8);
\draw (2,0) -- (2,8);
\draw (0,2) -- (4,2);
\draw (0,4) -- (4,4);
\draw (0,6) -- (4,6);
\node at (1,8.5) {\Large G};
\node at (3,8.5) {\Large D};
\node at (1,7) {11};
\node at (1,5) {13};
\node at (1,3) {17};
\node at (1,1) {19};
\node at (3,7) {21};
\node at (3,5) {11};
\node at (3,3) {19};
\node at (3,1) {9};
\node[anchor=east] at (7.2,7) {$11+21=$\textcolor{green!60!black}{32}};
\node[anchor=east] at (7.2,5) {$13+11=$\textcolor{green!60!black}{24}};
\node[anchor=east] at (7.2,3) {$17+19=$\textcolor{green!60!black}{36}};
\node[anchor=east] at (7.2,1) {$19+9=$\textcolor{green!60!black}{28}};
\node[anchor=east, magenta] at (7.2,-0.5) {*Pr4};
\draw (8,0) rectangle (10,8);
\draw (8,2) -- (10,2);
\draw (8,4) -- (10,4);
\draw (8,6) -- (10,6);
\node at (9,8.5) {\Large M};
\node at (9,7.3) {(4)};
\node at (9,6.7) {8};
\node at (9,5.3) {(4)};
\node at (9,4.7) {6};
\node at (9,3.3) {9};
\node at (9,2.7) {(4)};
\node at (9,1.3) {(4)};
\node at (9,0.7) {7};
\node[anchor=east] at (12.7,7) {$32+8=$\textcolor{green!60!black}{40}};
\node[anchor=east] at (12.7,5) {$24+6=$\textcolor{green!60!black}{30}};
\node[anchor=east] at (12.7,3) {$36+9=$\textcolor{green!60!black}{45}};
\node[anchor=east] at (12.7,1) {$28+7=$\textcolor{green!60!black}{35}};
\node[anchor=east, magenta] at (12.7,-0.5) {*Pr5};
\node at (1,-1.5) {60};   
\node at (2,-1.5) {$+$};  
\node at (3,-1.5) {60};   
\node at (6,-1.5) {$+$};  
\node at (9,-1.5) {30};   
\node at (11,-1.5) {$=$};
\node at (12,-1.5) {150};
\node[anchor=base west, magenta, inner sep=0pt]
  at (0,-2.5) {*Pr4:};
\node[anchor=base west, inner sep=0pt]
  at (1.1,-2.5)
  {$24 \;\rightarrow\; 28 \;\rightarrow\; 32
    \;\rightarrow\; 36$
    \qquad (common difference 4)};
\node[anchor=base west, magenta, inner sep=0pt]
  at (0,-3.4) {*Pr5:};
\node[anchor=base west, inner sep=0pt]
  at (1.1,-3.4)
  {$30 \;\rightarrow\; 35 \;\rightarrow\; 40
    \;\rightarrow\; 45$
    \qquad (common difference 5)};
\end{tikzpicture}
\caption{In the adjusted configuration, the G--D row sums form an arithmetic progression with common difference 4 when sorted. Adding the corresponding M primary values gives an arithmetic progression with common difference 5 when sorted. The bottom total, $60+60+30=150$, includes the G and D totals and only the four M primary values. It excludes the four secondary values of 4; the complete adjusted M column totals 46.}
\label{fig:Figure20}
\end{figure}

Furthermore, five main groupings spanning all three columns can be identified, as depicted in Figure~\ref{fig:Figure21}, each following the same rule: the outer pair (rows 1 and 4) and the inner pair (rows 2 and 3) yield the same sum. Within Columns G and D, two groupings per column each sum to 30 ($11+19$) and ($13+17$) in Column G, and ($21+9$) and ($11+19$) in Column D, while the cross-column groupings each sum to 60 ($32+28$ and $24+36$). In Column M, the two groupings of primary values each sum to 15 ($8+7$ and $6+9$). When the Column M primary values are added to their corresponding row sums, the resulting groupings each sum to 75 ($40+35$ and $30+45$).

A parallel pattern appears only within Column M when the subtraction and addition results are grouped according to the same rule, as shown in Figure~\ref{fig:Figure22}. Although the subtraction (each grouping summing to 7) and the addition side (each summing to 23) produce different individual sums, both sides follow the same grouping pattern, each pair of groupings summing to 30. When the primary values (6, 7, 8, and 9) of Column M are added to these grouped results, each group yields 45.

It was also found that when the doubled values (11 in Column G, 19 in Column D) and the secondary value 4 in Column M are set aside, the remaining unique numbers in each column appear to be distinguished by their mathematical properties. Column G contains two primes (13 and 17), and Column D contains two odd composites (9 and 21). Column M is mixed: 7 is prime, 9 is odd, and 6 and 8 are even.

\begin{figure}[H]
\centering
\begin{tikzpicture}[scale=0.80]
\useasboundingbox (-2,-1.5) rectangle (19,10.5);
\draw (0,0) rectangle (4,8);
\draw (2,0) -- (2,8);
\draw (0,2) -- (4,2);
\draw (0,4) -- (4,4);
\draw (0,6) -- (4,6);
\node at (1,8.5) {\Large G};
\node at (3,8.5) {\Large D};
\node at (1,7) {11};
\node at (1,5) {13};
\node at (1,3) {17};
\node at (1,1) {19};
\node at (3,7) {21};
\node at (3,5) {11};
\node at (3,3) {19};
\node at (3,1) {9};
\node at (1,-0.5) {60};
\node at (3,-0.5) {60};
\draw (0,7) -- (-1.2,7) -- (-1.2,1) -- (0,1);
\node[rotate=90] at (-1.5,4) {\small\textcolor{green!60!black}{30}};
\draw (0,5) -- (-0.55,5) -- (-0.55,3) -- (0,3);
\node[rotate=90] at (-0.85,4) {\small\textcolor{green!60!black}{30}};
\draw (4,7) -- (5.1,7) -- (5.1,1) -- (4,1);
\node[rotate=270] at (5.4,4) {\small\textcolor{green!60!black}{30}};
\draw (4,5) -- (4.5,5) -- (4.5,3) -- (4,3);
\node[rotate=270] at (4.8,4) {\small\textcolor{green!60!black}{30}};
\node at (6.8,7) {$11+21=$\textcolor{green!60!black}{$32$}};
\node at (6.8,5) {$13+11=$\textcolor{green!60!black}{$24$}};
\node at (6.8,3) {$17+19=$\textcolor{green!60!black}{$36$}};
\node at (6.8,1) {$19+9=$\textcolor{green!60!black}{$28$}};
\draw (8.7,7) -- (9.7,7) -- (9.7,1) -- (8.7,1);
\node[rotate=270] at (10.0,4) {\small\textcolor{green!60!black}{60}};
\draw (8.7,5) -- (9.1,5) -- (9.1,3) -- (8.7,3);
\node[rotate=270] at (9.4,4) {\small\textcolor{green!60!black}{60}};
\draw (10.8,0) rectangle (12.8,8);
\draw (10.8,2) -- (12.8,2);
\draw (10.8,4) -- (12.8,4);
\draw (10.8,6) -- (12.8,6);
\node at (11.8,8.5) {\Large M};
\node at (11.8,7.3) {(4)};
\node at (11.8,6.7) {\textcolor{green!60!black}{8}};
\node at (11.8,5.3) {(4)};
\node at (11.8,4.7) {\textcolor{green!60!black}{6}};
\node at (11.8,3.3) {\textcolor{green!60!black}{9}};
\node at (11.8,2.7) {(4)};
\node at (11.8,1.3) {(4)};
\node at (11.8,0.7) {\textcolor{green!60!black}{7}};
\draw (12.8,7) -- (13.8,7) -- (13.8,1) -- (12.8,1);
\node[rotate=270] at (14.1,4) {\small\textcolor{green!60!black}{15}};
\draw (12.8,5) -- (13.2,5) -- (13.2,3) -- (12.8,3);
\node[rotate=270] at (13.5,4) {\small\textcolor{green!60!black}{15}};
\node at (15.5,7) {$32+8=$\textcolor{green!60!black}{$40$}};
\node at (15.5,5) {$24+6=$\textcolor{green!60!black}{$30$}};
\node at (15.5,3) {$36+9=$\textcolor{green!60!black}{$45$}};
\node at (15.5,1) {$28+7=$\textcolor{green!60!black}{$35$}};
\draw (17.5,7) -- (18.5,7) -- (18.5,1) -- (17.5,1);
\node[rotate=270] at (18.8,4) {\small\textcolor{green!60!black}{75}};
\draw (17.5,5) -- (17.9,5) -- (17.9,3) -- (17.5,3);
\node[rotate=270] at (18.2,4) {\small\textcolor{green!60!black}{75}};
\end{tikzpicture}
\caption{Overview of the relationships among the groupings in the adjusted arrangement of Columns G, D, and M. The outer--inner grouping rule holds within Columns G and D and extends to their combined row sums: the outer pair ($32+28$) and the inner pair ($24+36$) each sum
to 60. Adding the corresponding primary values from Column M ($8,6,9,7$) to the row sums gives $40,30,45,35$. The resulting outer pair ($40+35$) and inner pair ($30+45$) each sum to 75.}
\label{fig:Figure21}
\end{figure}

\begin{figure}[H]
\centering
\begin{tikzpicture}[scale=0.80]
\useasboundingbox (-4.5,-1) rectangle (17,10.5);
\node at (-1.5,9) {Subtraction};
\node at (1,9) {\Large M};
\node at (3.5,9) {Addition};
\draw (0,0) rectangle (2,8);
\draw (0,2) -- (2,2);
\draw (0,4) -- (2,4);
\draw (0,6) -- (2,6);
\node[green!60!black] at (0.7,7) {8};
\node at (1.3,7) {4};
\node[green!60!black] at (0.7,5) {6};
\node at (1.3,5) {4};
\node[green!60!black] at (0.7,3) {9};
\node at (1.3,3) {4};
\node[green!60!black] at (0.7,1) {7};
\node at (1.3,1) {4};
\node at (-1.5,7) {$8-4=$\textcolor{green!60!black}{4}};
\node at (-1.5,5) {$6-4=$\textcolor{green!60!black}{2}};
\node at (-1.5,3) {$9-4=$\textcolor{green!60!black}{5}};
\node at (-1.5,1) {$7-4=$\textcolor{green!60!black}{3}};
\draw (-2.8,7) -- (-3.8,7) -- (-3.8,1) -- (-2.8,1);
\node[rotate=90] at (-4.1,4) {\small\textcolor{green!60!black}{7}};
\draw (-2.8,5) -- (-3.2,5) -- (-3.2,3) -- (-2.8,3);
\node[rotate=90] at (-3.5,4) {\small\textcolor{green!60!black}{7}};
\node at (3.5,7) {$8+4=$\textcolor{green!60!black}{12}};
\node at (3.5,5) {$6+4=$\textcolor{green!60!black}{10}};
\node at (3.5,3) {$9+4=$\textcolor{green!60!black}{13}};
\node at (3.5,1) {$7+4=$\textcolor{green!60!black}{11}};
\draw (5.2,7) -- (6.2,7) -- (6.2,1) -- (5.2,1);
\node[rotate=270] at (6.5,4) {\small\textcolor{green!60!black}{23}};
\draw (5.2,5) -- (5.6,5) -- (5.6,3) -- (5.2,3);
\node[rotate=270] at (5.9,4) {\small\textcolor{green!60!black}{23}};
\node at (8.0,7) {$4+12=$\textcolor{green!60!black}{16}};
\node at (8.0,5) {$2+10=$\textcolor{green!60!black}{12}};
\node at (8.0,3) {$5+13=$\textcolor{green!60!black}{18}};
\node at (8.0,1) {$3+11=$\textcolor{green!60!black}{14}};
\draw (9.8,7) -- (10.8,7) -- (10.8,1) -- (9.8,1);
\node[rotate=270] at (11.1,4) {\small\textcolor{green!60!black}{30}};
\draw (9.8,5) -- (10.2,5) -- (10.2,3) -- (9.8,3);
\node[rotate=270] at (10.5,4) {\small\textcolor{green!60!black}{30}};
\node at (12.5,7) {$16+8=$\textcolor{green!60!black}{24}};
\node at (12.5,5) {$12+6=$\textcolor{green!60!black}{18}};
\node at (12.5,3) {$18+9=$\textcolor{green!60!black}{27}};
\node at (12.5,1) {$14+7=$\textcolor{green!60!black}{21}};
\draw (14.3,7) -- (15.3,7) -- (15.3,1) -- (14.3,1);
\node[rotate=270] at (15.6,4) {\small\textcolor{green!60!black}{45}};
\draw (14.3,5) -- (14.7,5) -- (14.7,3) -- (14.3,3);
\node[rotate=270] at (15.0,4) {\small\textcolor{green!60!black}{45}};
\end{tikzpicture}
\caption{The grouping rule holds within Column M across both the subtraction and addition sides independently. On the subtraction side, both groupings sum to 7; on the addition side, both sum to 23. Although the two sides produce different totals, both follow the same pairing pattern, where the outer and inner pairs yield equal sums. When the subtraction and addition results are combined for each grouping, both yield 30; adding the corresponding primary value of Column M to each then yields 45.}
\label{fig:Figure22}
\end{figure}

The adjusted arrangement connects the three columns through complementary groupings and row-wise sum–product equalities. Section~\ref{sec:stateva} evaluates five specified structural properties under a shuffle model and across permitted readings and adjustments.

\section{Statistical Evaluation}
\label{sec:stateva}

The statistical evaluation was undertaken after the relationships and adjustments had been identified, so it serves to characterise the proposed structure, not to test a hypothesis specified in advance \parencite{popper1959}. It asks two questions. First, how often do shuffled arrangements of each reference configuration's values attain its composite score? This indicates how readily the defined patterns arise by chance when the values are placed at random. Second, which of the 64 permitted combinations of readings and adjustments satisfy the scored properties? This indicates how far the result depends on the specific choices made here. Both analyses are conditional on the observed values and are interpreted within the specified models.

\subsection{Composite Score and Structural Properties}
\label{subsec:stateva}

A composite score \(S\) was defined as the number of structural properties simultaneously satisfied by an arrangement:
\[
S=\sum_{j=1}^{5}
\mathbf{1}\!\left\{P_j\text{ is satisfied}\right\},
\]
where \(P_1,\ldots,P_5\) are the properties defined in Table~\ref{tab:properties}, and \(\mathbf{1}\{\cdot\}\) is the indicator function. Each property contributes one point, giving
a score between zero and five.

The score evaluates the operational definitions in Table~\ref{tab:properties}, rather than every relationship discussed in the exploratory analysis. The full set of arithmetic observations and the five scored properties should therefore be distinguished.

\begin{table}[H]
\centering
\normalsize
\renewcommand{\arraystretch}{1.15}
\setlength{\tabcolsep}{6pt}

\caption{The five structural properties evaluated in both the randomisation analysis and the sensitivity analysis of readings and adjustments. Each property is scored as present or absent.}
\label{tab:properties}

\begin{tabular}{@{}cp{12cm}@{}}
\toprule
\textbf{Property} &
\textbf{Operational definition} \\
\midrule
P1 &
The totals of Columns G and D are equal:
\(\sum_i G_i=\sum_i D_i\). \\

P2 &
The four row sums \(G_i+D_i\), when sorted, form an arithmetic progression with a nonzero common difference. \\

P3 &
For every row, the product of the two entries in the corresponding M pair equals the G--D row sum. \\

P4 &
Each of Columns G and D can be partitioned into two pairs with equal sums. Any pairing within each column is permitted; the pairs need not occupy the outer and inner positions. \\

P5 &
Taking the larger entry of each M pair as its primary value, the sum of the first and fourth primary values equals the sum of the second and third primary values. \\
\bottomrule
\end{tabular}
\end{table}

The score summarises simultaneous agreement with these five specified properties and does not assume that the properties are statistically independent. In the randomisation analysis, all five properties are evaluated together for every shuffled arrangement and only the composite score is reported, so any dependence among the properties is reflected in the simulated score distribution. In the enumeration of Section~\ref{sec:reconstruction-sensitivity}, each configuration is scored on the individual properties (Table~\ref{tab:reconstruction-property-results}) and on the composite score (Table~\ref{tab:reconstruction-score-distribution}). Because the properties and adjustments were identified through exploration of the same artefact, the shuffle analysis provides a fixed-score conditional benchmark, not a selection-adjusted test of the exploratory procedure \parencite{gelman2014statistical}.

\subsection{Reference Configurations}

The three reference configurations are defined in Table~\ref{tab:configurations}. Configurations 1 and 2 retain the documented positions and differ in their readings of the
damaged M pair. Configuration 3 is the fully adjusted arrangement developed in the preceding analysis.

\begin{table}[H]
\centering
\normalsize
\renewcommand{\arraystretch}{1.25}
\setlength{\tabcolsep}{6pt}

\caption{The three reference configurations evaluated in the randomisation analysis. Column entries and M pairs are listed in row order.}
\label{tab:configurations}

\begin{tabular}{@{}p{3.2cm}p{2.7cm}p{2.7cm}p{5cm}@{}}
\toprule
\textbf{Configuration} &
\textbf{Column G} &
\textbf{Column D} &
\textbf{Column M pairs} \\
\midrule
Config.\ 1\newline
(original, 9 and 4) &
\( (11,13,17,19) \) &
\( (11,21,19,9) \) &
\( (3,6),\ (4,8),\ (9,4),\ (5,7) \) \\

Config.\ 2\newline
(original, 10 and 5) &
\( (11,13,17,19) \) &
\( (11,21,19,9) \) &
\( (3,6),\ (4,8),\ (10,5),\ (5,7) \) \\

Config.\ 3\newline
(fully adjusted) &
\( (11,13,17,19) \) &
\( (21,11,19,9) \) &
\( (4,8),\ (4,6),\ (9,4),\ (4,7) \) \\
\bottomrule
\end{tabular}
\end{table}

Configurations 1 and 2 each satisfy P1 and P4 and therefore score \(S=2\). Configuration 3 satisfies all five properties and scores \(S=5\).

Although all three configurations contain sixteen entries, they do not contain the same multiset of numerical values. Configuration 2 replaces the damaged-pair reading \((9,4)\)
with \((10,5)\), while Configuration 3 replaces the separate entries \(3\) and \(5\) with two entries of \(4\). Each configuration consequently has its own randomisation distribution.

\subsection{Randomisation Model and Simulation}

For each reference configuration, its sixteen numerical entries were treated as a fixed multiset. In each simulated arrangement, the entries were shuffled uniformly into sixteen labelled positions: four positions in Column G, four in Column D, and eight positions forming four row-specific pairs in Column M. Repeated numerical values were retained with their original
multiplicities.

All five properties were evaluated using the resulting row correspondence. Reversing the two entries within an M pair does not change the score, but the row positions of the pairs
and the correspondence between Columns G, D, and M remain relevant.

For reference configuration \(c\), let \(H_{0,c}\) denote this shuffle model, and let \(S_{\mathrm{obs},c}\) be its observed score. The quantity of interest is the conditional tail
probability
\[
q_c=
\Pr_{H_{0,c}}
\left(S\geq S_{\mathrm{obs},c}\right).
\]
This is the probability that a shuffled arrangement achieves an equal or higher score under the specified model \parencite{Good2005, Ernst2004}.

A total of \(B=1{,}000{,}000\) shuffled arrangements were generated for each reference configuration. If \(K_c\) arrangements achieved a score at least as high as
\(S_{\mathrm{obs},c}\), the empirical tail frequency was
\[
\widehat{q}_c=\frac{K_c}{B}.
\]
The supplied Python implementation used a pseudorandom generator initialised with seed \(42\), processing the configurations in the order shown in Table~\ref{tab:configurations}.

The randomisation results reported here come from the simulations of one million arrangements. The additional runs of 50,000 arrangements produced by the accompanying script use the same sampling procedure, so they are not an independent validation.

The randomisation model is a conditional benchmark. It takes the values already present in each reference configuration as given and reshuffles only their positions \parencite{Fisher1935}. The proportions it reports therefore show how often shuffled arrangements of the same values produce at least as many patterns as the observed arrangement.

\subsection{Simulated Score Distribution}

Table~\ref{tab:scoredist} reports the score distribution obtained from the value pool of Configuration 3.

\begin{table}[H]
\centering
\normalsize
\renewcommand{\arraystretch}{1.15}
\setlength{\tabcolsep}{6pt}

\caption{Composite scores across 1,000,000 shuffled arrangements of Configuration 3's value pool. The Count and Proportion columns describe scores equal to \(s\); the final column
reports the empirical frequency of scores at least as high as \(s\).}
\label{tab:scoredist}

\begin{tabular}{@{}rrrr@{}}
\toprule
\textbf{Score \(s\)} &
\textbf{Count} &
\textbf{Proportion} &
\textbf{Tail frequency} \\
\midrule
0 & 919{,}486 & 91.9486\% & 100.0000\% \\
1 &  78{,}262 &  7.8262\% &   8.0514\% \\
2 &   2{,}202 &  0.2202\% &   0.2252\% \\
3 &      50  &  0.0050\% &   0.0050\% \\
4 &       0  &  0.0000\% &   0.0000\% \\
5 &       0  &  0.0000\% &   0.0000\% \\
\bottomrule
\end{tabular}
\end{table}

\subsection{Randomisation Results}

These results should be read as conditional numerical benchmarks. They show how unusual the observed scores are when the same numerical values are redistributed uniformly within the specified three-column structure.

\begin{table}[H]
\centering
\normalsize
\renewcommand{\arraystretch}{1.15}
\setlength{\tabcolsep}{6pt}

\caption{Randomisation results for the three reference configurations, using 1,000,000 shuffled arrangements per configuration. An exceedance is a simulated score at least as high as the observed score. Tail frequencies are conditional on each configuration's value pool and the
specified scoring rule.}
\label{tab:globalresults}

\begin{tabular}{@{}lcrr@{}}
\toprule
\textbf{Configuration} &
\textbf{Observed score} &
\textbf{Exceedances} &
\textbf{Tail frequency} \\
\midrule
Config.\ 1 (original, 9 and 4) &
2 & 2{,}805 & 0.2805\% \\

Config.\ 2 (original, 10 and 5) &
2 & 2{,}630 & 0.2630\% \\

Config.\ 3 (fully adjusted) &
5 & 0 & 0 observed \\
\bottomrule
\end{tabular}
\end{table}

For Configurations 1 and 2, approximately 0.26--0.28\% of shuffled arrangements achieved a score of at least two. These frequencies refer to $S\geq2$, not to scores of exactly two.

For Configuration 3, no shuffled arrangement achieved the observed score of five. Zero exceedances do not establish that the underlying tail probability is zero or below $1/B$. A conservative Monte Carlo estimate of the tail probability is $(K_c+1)/(B+1)$, which here equals $1/1{,}000{,}001\approx1.0\times10^{-6}$ \parencite{phipson2010permutation}.

Because the configurations have different value pools, these results are not a controlled comparison of positional adjustments alone. The following sensitivity analysis examines dependence on the specified readings and adjustments.

Under the binomial model for independent simulation trials, the one-sided exact 95\% upper confidence limit for a tail probability with zero exceedances \parencite{clopper1934use} is
\[
q_U=1-0.05^{1/B}.
\]
For \(B=1{,}000{,}000\), this gives
\[
q_U\approx2.996\times10^{-6},
\]
or approximately \(0.000300\%\). This limit quantifies simulation uncertainty under the fixed shuffle model and fixed scoring rule.

\subsection{Complete Enumeration of the Restricted Reading-and-Adjustment Space}
\label{sec:reconstruction-sensitivity}

A complete sensitivity analysis was conducted to determine how agreement with the structural properties depends on the specified readings and analytical adjustments. The documented
numerical arrangement and the ambiguity of the damaged entries in Column M provide the archaeological starting point \parencite[pp.~65--67]{deHeinzelin1957}. The positional exchanges and numerical adjustments are analytical choices proposed in the present study.

Following the general rationale of sensitivity analysis, this analysis examines the consequences of varying these inputs \parencite{SaltelliAnnoni2010}. Unlike the randomisation analysis, it preserves the column assignments and all values and positions except those explicitly permitted to vary. Its purpose is to identify which permitted combinations satisfy the properties.

\subsubsection{Definition of the Configuration Space}

Column G was fixed throughout as
\[
G=(11,13,17,19).
\]
Column D was permitted to take either its documented order or the order obtained by exchanging its first two entries:
\[
D\in
\left\{
(11,21,19,9),\,
(21,11,19,9)
\right\}.
\]

Column M was represented by four row-specific pairs. The first two pairs were permitted to occur in either order:
\[
\bigl((x,6),(4,8)\bigr)
\qquad\text{or}\qquad
\bigl((4,8),(x,6)\bigr),
\]
where \(x\in\{3,4\}\). The choice between \(3\) and \(4\) remained associated with the entry paired with \(6\), irrespective of the row occupied by that pair.

The third pair remained in the third row and was represented as
\[
(y,z),
\qquad
y\in\{9,10\},
\qquad
z\in\{4,5\}.
\]
All four combinations of these entries were included.

Their independent combination is an analytical choice for this sensitivity analysis.

The fourth pair remained in the fourth row and was represented as
\[
(w,7),
\qquad
w\in\{5,4\}.
\]

The configuration space therefore comprised six independently combinable binary choices, summarised in Table~\ref{tab:reconstruction-choices}. All other values and positions remained fixed.

\begin{table}[H]
\centering
\normalsize
\renewcommand{\arraystretch}{1.15}
\setlength{\tabcolsep}{6pt}

\caption{The six binary choices defining the restricted configuration space. Every combination was evaluated. Independent combinability is a feature of the enumeration.}
\label{tab:reconstruction-choices}

\begin{tabular}{@{}p{6cm}p{8cm}@{}}
\toprule
\textbf{Component} &
\textbf{Permitted alternatives} \\
\midrule
First two entries of D &
\( (11,21) \) or \( (21,11) \) \\

Order of the first two M pairs &
\( (x,6),(4,8) \) or \( (4,8),(x,6) \) \\

M entry paired with \(6\) &
\(x=3\) or \(x=4\) \\

First entry of the third M pair &
\(y=9\) or \(y=10\) \\

Second entry of the third M pair &
\(z=4\) or \(z=5\) \\

M entry paired with \(7\) &
\(w=5\) or \(w=4\) \\
\bottomrule
\end{tabular}
\end{table}

The number of distinct configurations was
\[
N=2^6=64.
\]

Every configuration was evaluated exactly once. No Monte Carlo sampling, random seed, or sampling-based confidence interval was required for this enumeration.

The value exchange in Column D and the exchange of the first two M pairs were allowed to occur separately or together. Likewise, the changes \(3\rightarrow4\) and \(5\rightarrow4\) were evaluated separately and together. The analysis therefore includes the coordinated single-unit reallocation described earlier but does not require that coordination or impose a fixed Column M total. It tests the consequences of retaining or omitting components of the proposed pathway, including configurations that do not preserve the total-preserving coordination of the unit transfer. These alternatives combine reading uncertainty with analytical adjustments.

\subsubsection{Evaluation and Results}

Each configuration was scored as in Subsection~\ref{subsec:stateva}, using the five properties defined in Table~\ref{tab:properties}. A complete match required \(S=5\). \\ 

Exactly one of the 64 configurations satisfied all five properties:
\[
\begin{aligned}
G &= (11,13,17,19),\\
D &= (21,11,19,9),\\
M &= \bigl((4,8),(4,6),(9,4),(4,7)\bigr).
\end{aligned}
\]
This is the fully adjusted reference configuration.

Individual property frequencies are reported in Table~\ref{tab:reconstruction-property-results}.

\begin{table}[H]
\centering
\normalsize
\renewcommand{\arraystretch}{1.15}
\setlength{\tabcolsep}{6pt}

\caption{Exact individual property frequencies across the 64 permitted configurations. Percentages describe the enumerated configuration space.}
\label{tab:reconstruction-property-results}

\begin{tabular}{@{}lrr@{}}
\toprule
\textbf{Property} &
\textbf{Passing configurations} &
\textbf{Percentage} \\
\midrule
P1 & \(64/64\) & 100.0000\% \\
P2 & \(32/64\) & 50.0000\% \\
P3 & \(1/64\)  & 1.5625\% \\
P4 & \(64/64\) & 100.0000\% \\
P5 & \(16/64\) & 25.0000\% \\
\midrule
All five simultaneously & \(1/64\) & 1.5625\% \\
\bottomrule
\end{tabular}
\end{table}

P1 and P4 held for every configuration because the permitted changes preserve the multisets of values in Columns G and D. Both columns retain their totals of \(60\) and their equal-sum pairings:
\[
11+19=13+17=30
\]
in Column G, and
\[
11+19=21+9=30
\]
in Column D. P1 and P4 are therefore invariant within this configuration space and do not distinguish between its alternatives.

P2 held in 32 configurations, corresponding to those in which Column D began with \(21\) followed by \(11\).
Their row sums were
\[
(32,24,36,28),
\]
which, when sorted, give
\[
(24,28,32,36),
\]
an arithmetic progression with common difference \(4\).
The documented D order instead gives
\[
(22,34,36,28),
\]
which does not form an arithmetic progression in any ordering.

P3 held in only one configuration. Its row-wise equalities were
\[
\begin{aligned}
11+21 &= 4\times8 = 32,\\
13+11 &= 4\times6 = 24,\\
17+19 &= 9\times4 = 36,\\
19+9  &= 4\times7 = 28.
\end{aligned}
\]
This was also the sole configuration satisfying all five properties. P5 held in 16 configurations, including the complete match, whose M primary values were \((8,6,9,7)\)
and satisfied

\[
8+7=6+9=15.
\]

The composite-score distribution is reported in Table~\ref{tab:reconstruction-score-distribution}. Both documented reference readings, with either \((9,4)\) or \((10,5)\) in the third M pair and without analytical adjustments, satisfied P1 and P4 only and therefore scored two out of five.

\begin{table}[H]
\centering
\normalsize
\renewcommand{\arraystretch}{1.15}
\setlength{\tabcolsep}{6pt}

\caption{Distribution of composite scores across the 64 permitted configurations. Each configuration is counted once.}
\label{tab:reconstruction-score-distribution}

\begin{tabular}{@{}crr@{}}
\toprule
\textbf{Score \(S\)} &
\textbf{Configurations} &
\textbf{Percentage} \\
\midrule
0 & 0  & 0.0000\% \\
1 & 0  & 0.0000\% \\
2 & 24 & 37.5000\% \\
3 & 32 & 50.0000\% \\
4 & 7  & 10.9375\% \\
5 & 1  & 1.5625\% \\
\midrule
Total & 64 & 100.0000\% \\
\bottomrule
\end{tabular}
\end{table}

The supplementary reconstruction script exports the results as CSV tables, listing all 64 configurations, their six choices, the five property outcomes, and the composite scores. It also reports how many configurations satisfy each pair of properties and how each of the six choices affects the property outcomes (Subsection~\nameref{subsec:DandCA}).

\subsubsection{Scope and Limitations}
\label{sec:scopeandlimitations}

The fraction \(1/64\) is an exact descriptive proportion within the specified configuration space, not a randomisation tail probability or a historical probability of chance or intention. Its value depends on the alternatives included in the enumeration.

The complete match does not represent five independent restrictions converging on the same arrangement. P1 and P4 hold throughout the permitted space, and P3 alone selects the configuration satisfying all five properties. The counts therefore describe overlapping structural relationships, not five independent pieces of evidence.

The configuration space permits only the documented order in Column~D and the exchange of its first two entries. Exchanging the last two entries instead gives \(D=(11,21,9,19)\), which also yields equal outer and inner sums (\(11+19=21+9=30\)). Its G--D row sums are \((22,34,26,38)\), which do not form an arithmetic progression, so P2 is not satisfied. This alternative is
not part of the 64-configuration enumeration, and the results do not establish uniqueness over all possible exchanges in Column~D. Evaluating the alternative exchanges is left for further research.

\subsection{Interpretation}

The two analyses answer different questions. The randomisation analysis asks how readily the defined patterns arise when a configuration's own values are shuffled; the enumeration asks how far the result depends on the readings and adjustments made here. For Configurations~1 and~2, the score of two reflects only the equal totals and pairings of Columns~G and~D (P1 and P4), not any adjustment. For Configuration~3, no shuffled arrangement reached the observed score
(Table~\ref{tab:globalresults}), and in the enumeration the row-wise product property (P3) alone selects the single complete match. Neither analysis tests whether the arrangement was produced by intention, and, because the properties were identified after the adjustments, neither is a prospective test \parencite{diaconis1989, gelman2014statistical}.

Subsection~\nameref{subsec:methodological-limitations}

\section{Storytelling with Numbers}
\label{sec:storytelling}

The arithmetic relationships identified above motivate consideration of how such an arrangement, if it were intentional, could have been used to demonstrate numerical ideas. The following scenario is motivated by the demonstrable mathematical relationships and the logic of the adjustments.

One possibility is that the notch groups served as a reference for arranging physical markers. In such a setting, positional exchanges and a unit transfer could make numerical equivalences visible and permit a sequence of demonstrations. Figure~\ref{fig:Figure23} illustrates
this interpretation. The proposed flexibility belongs to the external arrangement of markers, inspired by the adjustments made and the hypothesis that the notches represent flexible values that needed to be laid out.

\begin{figure}[H]
\centering
\begingroup
\definecolor{cupoutline}{RGB}{105,91,74}
\definecolor{cupfill}{RGB}{247,239,224}
\definecolor{marblefill}{RGB}{112,153,169}
\definecolor{marbleoutline}{RGB}{65,103,120}
%
\newcommand{\schematiccup}[4]{%
    \pgfmathsetmacro{\cupwidth}{#1}
    \pgfmathsetmacro{\cupheight}{#2}
    \pgfmathsetmacro{\marbleradius}{0.055*\cupwidth}
    \pgfmathsetmacro{\marblestep}{0.66*\cupwidth/#4}
    \path[
        fill=cupfill,
        draw=cupoutline,
        line width=0.8pt,
        rounded corners=2pt
    ]
        ({-0.50*\cupwidth},{0.50*\cupheight})
        --
        ({-0.36*\cupwidth},{-0.50*\cupheight})
        --
        ({0.36*\cupwidth},{-0.50*\cupheight})
        --
        ({0.50*\cupwidth},{0.50*\cupheight})
        -- cycle;
    \foreach \n in {1,...,#3}{%
        \pgfmathtruncatemacro{\marblerow}
            {floor((\n-1)/#4)}
        \pgfmathtruncatemacro{\marblecol}
            {mod(\n-1,#4)}
        \pgfmathsetmacro{\marblex}
            {(\marblecol-(#4-1)/2)*\marblestep}
        \pgfmathsetmacro{\marbley}
            {-0.50*\cupheight
             +1.65*\marbleradius
             +2.25*\marbleradius*\marblerow}
        \filldraw[
            fill=marblefill,
            draw=marbleoutline,
            line width=0.3pt
        ]
            (\marblex,\marbley)
            circle[radius=\marbleradius];
    }
    \draw[
        cupoutline,
        line width=0.8pt,
        fill=cupfill!35!white
    ]
        (0,{0.50*\cupheight})
        ellipse[
            x radius={0.50*\cupwidth},
            y radius={0.075*\cupheight}
        ];
}
\begin{tikzpicture}[
    x=1cm,
    y=1cm,
    every node/.style={font=\small},
    quantity/.style={
        anchor=west,
        font=\small,
        text=black
    }
]
\node[font=\large] at (0,1.30) {G};
\node[font=\large] at (3.15,1.30) {D};
\node[font=\large] at (6.15,1.30) {M};

\foreach \y in {0,-2.25,-4.50,-6.75}{
    \draw[
        gray!35,
        densely dashed,
        line width=0.4pt
    ]
        (-1.05,{\y-0.98})
        --
        (7.35,{\y-0.98});
}
\foreach \quantity/\y in {
    11/0,
    13/-2.25,
    17/-4.50,
    19/-6.75
}{
    \begin{scope}[shift={(0,\y)}]
        \schematiccup{1.65}{1.05}{\quantity}{7}
        \node[quantity] at (1.00,0) {\quantity};
    \end{scope}
}
\foreach \quantity/\y in {
    21/0,
    11/-2.25,
    19/-4.50,
    9/-6.75
}{
    \begin{scope}[shift={(3.15,\y)}]
        \schematiccup{1.65}{1.05}{\quantity}{7}
        \node[quantity] at (1.00,0) {\quantity};
    \end{scope}
}
\foreach \upperquantity/\lowerquantity/\y in {
    4/8/0,
    4/6/-2.25,
    9/4/-4.50,
    4/7/-6.75
}{
    \begin{scope}[shift={(6.15,{\y+0.48})}]
        \schematiccup{1.05}{0.62}{\upperquantity}{3}
        \node[quantity] at (0.68,0) {\upperquantity};
    \end{scope}

    \begin{scope}[shift={(6.15,{\y-0.48})}]
        \schematiccup{1.05}{0.62}{\lowerquantity}{3}
        \node[quantity] at (0.68,0) {\lowerquantity};
    \end{scope}
}
\end{tikzpicture}
\endgroup
\caption{Hypothetical representation of the adjusted numerical configuration using cups and marbles. Columns G and D contain one larger cup per row, while Column M contains two smaller cups per corresponding row. Each cup represents one notch group, and each marble represents one unit. The labels indicate the quantities represented. This illustrates how the notch groups of the Ishango Bone might have served as a reference for a dynamic physical arrangement. Such a setup would have allowed quantities to be adjusted during a mathematical demonstration. }
\label{fig:Figure23}
\end{figure}

Physical markers, such as stones arranged on leather sheets or in cups, may have represented numerical values in a dynamic, interactive way, allowing quantities to be adjusted during demonstrations. Such an idea is compatible with proposals that manipulable material forms scaffold number concepts \parencite{overmann2013material}, and with the archaeological record
of counters and tokens in later periods \parencite{schmandtbesserat1992before}.

\subsection*{Possible Storytelling Purpose and Sequence}

If the arithmetic relationships identified here were intentional, one possible setting is a shared mathematical practice, and the Ishango Bone may then represent more than an individual's record of quantities. It could bring together numerical ideas explored, shared, and refined through interaction, perhaps across generations. Although the operations themselves are elementary, their interweaving within a small arrangement gives the proposed structure a
degree of economy that invites consideration of such a practice.

Jean de Heinzelin concluded from his study of the Ishango area that it represented a relatively isolated settlement. The exceptional conditions at Ishango provided an abundant and reliable food supply, creating a secure environment that enabled a sedentary way of life
\parencite[pp.~91--92]{deHeinzelin1957}. Such conditions may have offered a sheltered setting in which numerical ideas could be developed. This also raises the question of the setting in which a mathematical concept of this kind might have been used.

On this interpretation, the artefact would not necessarily mark the invention of the relationships it expresses. It could instead reflect shared knowledge of elementary arithmetic relationships and a tradition of passing such knowledge on. This section develops one possible way in which the bone might have been used, taking the adjusted arrangement described above as its starting point.

\begin{figure}[H]
\centering
\begin{tikzpicture}[scale=1]
\useasboundingbox (-0.5,-1) rectangle (10,9.5);
\draw (0,0) rectangle (4,8);
\draw (2,0) -- (2,8);
\draw (0,2) -- (4,2);
\draw (0,4) -- (4,4);
\draw (0,6) -- (4,6);
\node at (1,8.5) {\Large G};
\node at (3,8.5) {\Large D};
\node at (1,7) {11};
\node at (1,5) {13};
\node at (1,3) {17};
\node at (1,1) {19};
\node at (3,7) {11};
\node at (3,5) {21};
\node at (3,3) {19};
\node at (3,1) {9};
\draw (6,0) rectangle (8,8);
\draw (6,2) -- (8,2);
\draw (6,4) -- (8,4);
\draw (6,6) -- (8,6);
\node at (7,8.5) {\Large M};
\node at (7,7.3) {3};
\node at (7,6.7) {6};
\node at (7,5.3) {4};
\node at (7,4.7) {8};
\node at (7,3.3) {9};
\node at (7,2.7) {4};
\node at (7,1.3) {5};
\node at (7,0.7) {7};
\end{tikzpicture}
\caption{Numerical entries of Columns G, D, and M in their documented order, using the provisional reading $(9,4)$ for the damaged pair in Column M. The M entries are displayed in the analytical pairs used in this study. This configuration serves as the starting point for the hypothetical demonstration described below.}
\label{fig:Figure24}
\end{figure}

The presenter might have started with the secondary values, numbers 4, in Column~M, as shown in Figure~\ref{fig:Figure25}, demonstrating to the audience that different combinations of numbers can yield the same total. In teaching terms, this is a concrete demonstration of numerical equivalence and of the composition of a number ($8=3+5=4+4$). This could also have served as an early introduction to the grouping principle, the pairing of outer and inner values, that recurs
across the entire arrangement.

\begin{figure}[H]
\centering
\begin{tikzpicture}[scale=1]
\useasboundingbox (-0.5,-1) rectangle (4.5,9.5);
\draw (0,0) rectangle (2,8);
\draw (0,2) -- (2,2);
\draw (0,4) -- (2,4);
\draw (0,6) -- (2,6);
\node at (1,8.5) {\Large M};
\node[green!60!black] at (1,7.3) {3};
\node at (1,6.7) {(6)};
\node[green!60!black] at (1,5.3) {4};
\node at (1,4.7) {(8)};
\node at (1,3.3) {(9)};
\node[green!60!black] at (1,2.7) {4};
\node[green!60!black] at (1,1.3) {5};
\node at (1,0.7) {(7)};
\draw (2,7) -- (3.0,7) -- (3.0,1) -- (2,1);
\node[green!60!black][rotate=270] at (3.3,4) {8};
\draw (2,5) -- (2.5,5) -- (2.5,3) -- (2,3);
\node[green!60!black][rotate=270] at (2.75,4) {8};
\end{tikzpicture}
\caption{Using the original entries in Column M, the secondary values (3, 4, 4, and 5) grouped by the overall grouping pattern appear to produce the same total: the outer pair (3 + 5) and the inner pair (4 + 4) both equal 8. If one unit is taken from 5 and added to 3, all secondary numbers become 4.}
\label{fig:Figure25}
\end{figure}

In this way, an audience could observe directly that 3 + 5 equals 4 + 4. Moving a single stone from the group of 5 to the group of 3 would make the equivalence visible and show that all four secondary values become 4. Such a demonstration might have conveyed a conceptual understanding of numerical equivalence, a theme that recurs across the arrangement: different sums yield the same result, and the products of the M pairs equal the row sums of Columns~G and~D. The quartz-tipped end of the bone could have served as a pointer or stylus, for example, for drawing grouping connections in sand.

It remains an open question why the documented positions in Column~D and the corresponding two pairs in Column~M differ from the arrangement in which the cross-column relationships appear. Once these repositioning steps had been made, however, the relationships across all three columns could have been demonstrated in many different sequences, in whatever order suited the audience or context. All adjustments are illustrated in Figure~\ref{fig:Figure26}.

\begin{figure}[H]
\centering
\begin{tikzpicture}[scale=1]
\useasboundingbox (-0.5,-1) rectangle (10.5,9.5);
\draw (0,0) rectangle (4,8);
\draw (2,0) -- (2,8);
\draw (0,2) -- (4,2);
\draw (0,4) -- (4,4);
\draw (0,6) -- (4,6);
\node at (1,8.5) {\Large G};
\node at (3,8.5) {\Large D};
\node at (1,7) {11};
\node at (1,5) {13};
\node at (1,3) {17};
\node at (1,1) {19};
\node at (3,7) {11};
\node at (3,5) {21};
\node at (3,3) {19};
\node at (3,1) {9};
\draw (6,0) rectangle (8,8);
\draw (6,2) -- (8,2);
\draw (6,4) -- (8,4);
\draw (6,6) -- (8,6);
\node at (7,8.5) {\Large M};
\node at (7,7.3) {3\textcolor{green!60!black}{(+1)}};
\node at (7,6.7) {6};
\node at (7,5.3) {4};
\node at (7,4.7) {8};
\node at (7,3.3) {9\textcolor{black}{(+1?)}};
\node at (7,2.7) {4\textcolor{black}{(+1?)}};
\node at (7,1.3) {5\textcolor{green!60!black}{(-1)}};
\node at (7,0.7) {7};
\draw[green!60!black] (4.0,7) -- (4.7,7) -- (4.7,5) -- (4.0,5);
\fill[green!60!black] (4.0,7) -- (4.15,7.08) -- (4.15,6.92);
\fill[green!60!black] (4.0,5) -- (4.15,5.08) -- (4.15,4.92);
\draw[green!60!black] (8.0,7) -- (8.7,7) -- (8.7,5) -- (8.0,5);
\fill[green!60!black] (8.0,7) -- (8.15,7.08) -- (8.15,6.92);
\fill[green!60!black] (8.0,5) -- (8.15,5.08) -- (8.15,4.92);
\draw[green!60!black, dashed] (8.0,7.3) -- (9.0,7.3) -- (9.0,1.3) -- (8.0,1.3);
\fill[green!60!black] (8.0,7.3) -- (8.15,7.38) -- (8.15,7.22);
\fill[green!60!black] (8.0,1.3) -- (8.15,1.38) -- (8.15,1.22);
\node[rotate=270] at (9.3,4.3) {\small\textcolor{green!60!black}{$\pm$1}};
\end{tikzpicture}
\caption{An overview of all exploratory adjustments implemented in this study, suggesting how these adjustments might have been introduced within the storytelling sequence. First, one unit is reallocated from 5 to 3 in Column M (dashed bracket); this is then followed by the positional exchange of numbers 21 and 11 in Column D and the corresponding exchange of number pairs (4, 6) and (4, 8) in Column M (solid brackets).}
\label{fig:Figure26}
\end{figure}

Viewed in this light, the deviations from the documented arrangement (the positional exchange in Column~D, the exchange of two pairs in Column~M, and the unit transfer) need not be corrections of transmission errors, as first assumed. They could instead have been steps in a demonstration carried out with physical markers. This offers an alternative explanation of the adjustments.

\section{Discussion}
\label{sec:discussion}

The analysis distinguishes properties of the documented transcription from relationships obtained through reconstruction. The former include equal G and D totals, complementary pairings to 30, and the representation of every odd integer from 9 to 21 except 15. The latter connect all three columns through the proposed positional exchanges and unit transfer. The question is whether this extension reflects a deliberate organisation or analytical choice that produced it.

It seems also likely that Column~M is a later addition, as it relates to and complements Columns~G and~D, which already form a noteworthy arrangement in their own right. Further traceological (microwear) analyses could clarify whether all three columns were produced at the same time or whether Columns~G and~D reflect an earlier stage in the development of the arrangement.

\subsection*{Addressing Keller's Critique of Mathematical Intent in Tally Bones}

Keller cautions against interpreting mathematical relationships identified in prehistoric markings as evidence of the makers' mathematical intentions \parencite{Keller2010}. His central objection is not that such relationships cannot be found, but that assigning numerical meanings to the marks does not establish which, if any, of them reflects their original purpose. This
distinction between a mathematical description and a historical interpretation is relevant to the present study.

Although some of the numerical features examined here have been mentioned in previous studies, they have not yet been subjected to systematic documentation or detailed examination. De Heinzelin himself drew attention to internal arithmetic relationships in the documented arrangement and wrote:

\begin{quote}
Furthermore, a series of internal relationships make this table a fascinating puzzle where you can never be sure you’ve exhausted all the possible combinations: duplicate numbers, products equal to sums, sums equal to prime numbers, sums equal to the 4 times table, column sums equal to 60, and so on. \parencite[p.~67]{deHeinzelin1957}
\end{quote}

\noindent
The English rendering is a machine translation of the French original, produced with DeepL Translator \parencite{DeepL2026} on 26 November 2026.\footnote{Original French: \enquote{En outre, une série de relations internes font de ce tableau un jeu passionnant dont on n'est jamais sûr d'avoir épuisé toutes les combinaisons : duplication des nombres, produits égaux à des sommes, sommes égales à des nombres premiers, sommes égales à la table de 4, additions de colonnes égales à 60, etc.} \parencite[p.~67]{deHeinzelin1957}.} \\ 

Keller later described the same complementary pairings summing to 30 in Columns~G and~D, which were also identified independently in the present study:

\begin{quote}
And here's another potential \enquote{game}: on the left, the two end numbers total 30, as do the two middle numbers, whereas on the right the first and the third numbers total 30, and likewise the second and the fourth. \parencite[p.~5]{Keller2010}
\end{quote}

Keller, however, rejects such pairings and mathematical relationships as sufficient evidence of intention, while pointing to ethnographic cases where grouped marks admit a mathematical description although they were not meant as arithmetic demonstrations \parencite[pp.~13--16]{Keller2010}. Nevertheless, the possibility that numerical patterns can arise incidentally does not, by itself, justify dismissing the particular relationships examined here.

Questions of intent are unlikely to be settled by the numerical pattern alone and would need to be considered together with material, contextual and comparative evidence. The explicit rationale for the adjustments, together with the cross-column arithmetic relationships and grouping rules that connect the adjusted arrangement, nevertheless offers reasonable grounds to retain intentionality as a hypothesis worthy of further investigation.

\subsection*{Open Questions and Directions for Further Research}

Two structural features of the arrangement remain unexplained and warrant further investigation, should the broader arrangement indeed reflect an intentional design.

First, it can be observed that reordering the fully adjusted arrangement by placing rows in the sequence 2, 4, 1, 3 rather than 1, 2, 3, 4, as shown in Figure~\ref{fig:Figure27}, preserves all grouping rules and multiplicative, as well as division, relationships across all three columns. This alternative ordering would additionally produce seamless sequential progressions throughout: the row sums of Columns G and D would appear in consecutive order from top to bottom ($24 \rightarrow 28 \rightarrow 32 \rightarrow 36$), as well as the primary numbers in Column M; the arithmetic progressions \textit{Pr4} and \textit{Pr5} would align sequentially; and the addition and subtraction results within Column M would likewise appear in consecutive order. This raises a question: if the arrangement was intentional, why was it ordered this way rather than in fully sequential order? One possibility is that the order was chosen to preserve some other functionality or pattern that might not yet have been identified in this analysis. Another is that it might reflect an aspect of the artefact's practical use that remains to be understood. A further possibility is pedagogical: a non-sequential order may have been chosen deliberately to introduce a degree of difficulty, so that an audience could not easily anticipate the results during teaching.

\begin{figure}[H]
\centering
\begin{tikzpicture}[scale=1, every node/.style={font=\large}]
\useasboundingbox (-1.7,-1.8) rectangle (14.2,9.5);
\begin{scope}
  \draw (0,0) rectangle (6,8);
  \draw (2,0) -- (2,8);
  \draw (4,0) -- (4,8);
  \draw (0,2) -- (6,2);
  \draw (0,4) -- (6,4);
  \draw (0,6) -- (6,6);
  \node at (1,8.5) {\LARGE G};
  \node at (3,8.5) {\LARGE D};
  \node at (5,8.5) {\LARGE M};
  \node at (-1,7) {Row 1};
  \node at (-1,5) {Row 2};
  \node at (-1,3) {Row 3};
  \node at (-1,1) {Row 4};
  \node at (1,7) {11};
  \node at (1,5) {13};
  \node at (1,3) {17};
  \node at (1,1) {19};
  \node at (3,7) {21};
  \node at (3,5) {11};
  \node at (3,3) {19};
  \node at (3,1) {9};
  \node at (5,7.25) {4};
  \node at (5,6.75) {8};
  \node at (5,5.25) {4};
  \node at (5,4.75) {6};
  \node at (5,3.25) {9};
  \node at (5,2.75) {4};
  \node at (5,1.25) {4};
  \node at (5,0.75) {7};
  \node at (3,-1) {(a)};
\end{scope}
\begin{scope}[xshift=8cm]
  \draw (0,0) rectangle (6,8);
  \draw (2,0) -- (2,8);
  \draw (4,0) -- (4,8);
  \draw (0,2) -- (6,2);
  \draw (0,4) -- (6,4);
  \draw (0,6) -- (6,6);
  \node at (1,8.5) {\LARGE G};
  \node at (3,8.5) {\LARGE D};
  \node at (5,8.5) {\LARGE M};
  \node at (-1,7) {Row 2};
  \node at (-1,5) {Row 4};
  \node at (-1,3) {Row 1};
  \node at (-1,1) {Row 3};
  \node at (1,7) {13};
  \node at (1,5) {19};
  \node at (1,3) {11};
  \node at (1,1) {17};
  \node at (3,7) {11};
  \node at (3,5) {9};
  \node at (3,3) {21};
  \node at (3,1) {19};
  \node at (5,7.25) {4};
  \node at (5,6.75) {6};
  \node at (5,5.25) {4};
  \node at (5,4.75) {7};
  \node at (5,3.25) {4};
  \node at (5,2.75) {8};
  \node at (5,1.25) {9};
  \node at (5,0.75) {4};
  \node at (3,-1) {(b)};
\end{scope}
\end{tikzpicture}
\caption{Values of Columns G, D, and M. (a) Original row order used in this study. (b) Rows reordered as 2, 4, 1, 3; all grouping rules and operations are retained, while the numerical progressions, including addition and subtraction sequences, would appear in consecutive order.}
\label{fig:Figure27}
\end{figure}

Second, within the adjusted Column M arrangement, the third-row pair, provisionally read as $(9,4)$, is the only pair in which the primary value precedes the secondary value 4. The other three pairs place the secondary value first. If the primary/secondary ordering elsewhere in Column M is intentional, this single reversal would constitute a genuine anomaly requiring explanation. Experimental archaeology may help assess their practical implications within the proposed demonstrative arrangement.

\subsection*{Intentional and Incidental Patterns}

The aim of this study was to identify patterns, and the deviations from them, that involve every one of the sixteen numerical entries across the three columns. Earlier work also proposed arithmetic connections spanning all three columns, including the calculation-rule interpretation of Pletser and Huylebrouck \parencite{Pletser1999}. The present analysis examines a different proposed set of cross-column relationships, developed through the positional exchanges and unit transfer described above.

The proposed arrangement brings together elementary arithmetic operations through an interwoven set of cross-column relationships and grouping rules. In the adjusted configuration, the outer--inner grouping rule produces equal sums within Columns G and D and within the primary values of Column M, as well as arithmetic operations and their results. It should also be noted that if interchanging Columns G and D, these numerical relationships remain preserved. After the Column D exchange, the G--D row sums are $(32,24,36,28)$ from top to bottom. Arranged in consecutive order, they become $(24,28,32,36)$, an arithmetic progression with common difference 4. Adding the corresponding Column M primary values to those row sums gives a second arithmetic progression, with common difference 5 when sorted. The primary values of Column M, together with the results of adding and subtracting the secondary value 4, span the consecutive integers from 2 to 13. The operations are elementary, but their interweaving gives the proposed arrangement a coherent mathematical structure.

Several further features suggest that the number 4 itself may have served as an organising principle, and not only emerged as a byproduct of the individual column patterns. For example, Columns~G and~D each contain four numbers, and Column~M contains eight, giving sixteen numbers in total. Columns~G and~D each sum to 60, a multiple of 4 ($60=4\times15$), and their combined total of 120 is likewise a multiple of 4 ($120=4\times30$). Whether this recurrence of 4 reflects a deliberate compositional choice or is an arithmetic consequence of the numbers selected remains unclear.

The analysis also identified further patterns whose intentionality is harder to assess. These include the diagonal placement of the doubled primes 11 and 19 (Figure~\ref{fig:Figure8}), the equal sums obtained by adding each unique number to the doubled number in its row and in its column (Figure~\ref{fig:Figure9}), and the equal block totals of 45 obtained by subtracting one occurrence of the doubled number from each block (Figure~\ref{fig:Figure11}). The last of these relationships already holds in the documented arrangement. In the adjusted arrangement, the same four values ($24,28,32,36$) can be obtained in two ways, as shown in Figure~\ref{fig:Figure28}: as the sums of the four rows of Columns~G and~D, and as the sums of the upper and lower pairs of entries within each column ($11+13=24$ and $17+19=36$ in Column~G; $21+11=32$ and $19+9=28$ in Column~D). All four values are multiples of 4; consequently, the differences between them ($4,8,12$) are also multiples of 4. Although these patterns involve all entries in Columns~G and~D, it remains unclear whether they were intended or arose as byproducts of the selected numbers and their arrangement.

\begin{figure}[H]
\centering
\begin{tikzpicture}[scale=0.80]
\useasboundingbox (-0.5,-1) rectangle (10,9.5);
\draw (0,5.5) rectangle (3,8.5);
\node at (1.5,7) {\large 24};
\draw (6,5.5) rectangle (9,8.5);
\node at (7.5,7) {\large 32};
\draw (0,0.5) rectangle (3,3.5);
\node at (1.5,2) {\large 36};
\draw (6,0.5) rectangle (9,3.5);
\node at (7.5,2) {\large 28};
\node at (1.5,9.2) {\Large G};
\node at (7.5,9.2) {\Large D};
\coordinate (C) at (4.5,4.5);
\draw (3,7) -- node[above=6pt]{\large 8} (6,7);
\draw (3,2) -- node[below=6pt]{\large 8} (6,2);
\draw (1.5,5.5) -- node[left=6pt]{\large 12} (1.5,3.5);
\draw (7.5,5.5) -- node[right=6pt]{\large 4} (7.5,3.5);
\draw (3,5.5) -- (C) -- (6,3.5);
\draw (6,5.5) -- (C) -- (3,3.5);
\node[above=6pt] at (C) {\large 4};
\node[below=6pt] at (C) {\large 4};
\end{tikzpicture}
\caption{In Columns G and D, the four values produced by pairing each unique number with its adjacent doubled number (24, 28, 32, 36), as previously shown in Figure~\ref{fig:Figure9}, produce a further pattern when arranged vertically and as four separate blocks. Not only are the values themselves all multiples of 4, but the differences along all connecting lines, 4, 8, and 12, are likewise all multiples of 4.}
\label{fig:Figure28}
\end{figure}

Additional patterns may remain to be identified within this arrangement. Experimental archaeology could assess whether the proposed demonstrations are practically feasible and how the patterns might have been recognised and used. Such work could show that the relationships were workable, but not that the makers knew of or intended them, since it is unlikely that all patterns have been chosen intentionally and some may be incidental byproducts of the numbers chosen \parencite{diaconis1989}.

\subsection*{Mathematical Knowledge and Oral Transmission}

Were this numerical arrangement deliberate, it could have supported more than the demonstration of individual calculations. Different quantities and operations could be visually demonstrated to produce equivalent results or preserve the same grouping relationships, making numerical equivalence and structure available for exploration. On this interpretation, the proposed
positional exchanges and unit transfer would correspond to steps in such a demonstration rather than to corrections of the artefact, and the bone might have served as a reference for active demonstration, rather than solely as a fixed record of quantities.

Such a use would be consistent with a social practice of mathematical communication, embedded in its cultural setting \parencite{dambrosio1985}. If the bone served as a reference for teaching, it could have been used to show numerical ideas that were already communicated through speech, gesture, and arrangements of objects \parencite{kelly2015}. Whether such a practice relates to later mathematical traditions is a separate question.

No documented line of continuity connects the Ishango Bone to later written mathematics. The earliest surviving written numerical records, the proto-cuneiform accounting tablets from Uruk (late fourth millennium BCE), already attest to well-established practices of accounting and measurement \parencite{Englund2004}, and Old Babylonian texts contain procedures equivalent
to solving quadratic equations \parencite{Joseph2011, Neugebauer1969}. These records are separated from Ishango by well over 15{,}000 years, as well as by the geographic distance between central Africa and Mesopotamia. Moreover, they document practices that were already in use rather than the first emergence of elementary arithmetic reasoning. Within Africa, notched bones such as the Lebombo bone from southern Africa are older than Ishango, and ethnographic
accounts describe the use of counting sticks and strings in central Africa \parencite{Huylebrouck2006}. However, no comparable artefacts have to date been documented from the intervening periods in this region. Whether any link exists therefore remains an open question, and the possibility that relevant artefacts have been lost, overlooked, or not yet discovered cannot be excluded.

Against this background, the possibility that the Ishango Bone is a much earlier expression of arithmetic reasoning cannot be excluded. Its potential significance lies in the apparent combination of arithmetic operations, complementary groupings, and numerical regularities within a single artefact. If this organisation was deliberate, the bone might preserve an early example of mathematical concepts arranged for demonstration and possibly shared within a community. It would not necessarily represent the origin of those concepts, nor a necessary stage in a single line of development leading to later mathematics. It could be an isolated, regionally and temporally limited development that was later rediscovered independently, or it could reflect one of the earliest stages of mathematical reasoning in humans more generally. Without further comparable artefacts, the evidence does not allow a choice between these possibilities.

The gaps in the surviving record leave the subsequent history of any such practice unresolved. If numerical knowledge was transmitted primarily through oral tradition, gesture, and demonstrations using perishable objects such as bones or wood, which rarely survive intact, such practices may have left few recognisable archaeological traces \parencite{kelly2015}. This possibility leaves room for continuity without establishing it. The Ishango Bone nevertheless remains of interest in its own context, as a possible early material reference for exploring arithmetic relationships and communicating mathematical ideas. It also highlights the need for more reliable frameworks for assessing intentionality in marked artefacts.

\subsection*{Methodological Scope and Limitations}
\label{subsec:methodological-limitations}

The analysis draws on de Heinzelin’s transcription, a provisional reading of the damaged entries, a proposed correspondence between the number pairs in Column M and the rows of Columns G and D, and explicit numerical and positional adjustments. Stating these choices explicitly allows the proposed structure, and the dependence of its patterns on them, to be examined.

The quantitative analyses describe the numerical structure of the proposed arrangement and the conditions under which it holds. The shuffle model compares the arrangement only with random placement of its values. The enumeration varies only the readings and adjustments proposed here, so it shows how the result depends on those choices. Neither analysis can therefore settle whether the arrangement was intentional.

\section{Conclusion}
\label{sec:conclusion}

This study documents a connected arithmetic analysis of the Ishango Bone and develops a possible demonstrative use. The documented arrangement already displays equal totals of 60 in Columns~G and~D and subgroup-pairings summing to 30. The proposed positional exchanges in Columns~D and~M, together with a provisional reading of the damaged pair and a total-preserving unit transfer, define an adjusted arrangement that connects all three columns through a shared outer--inner grouping rule, an arithmetic progression of row sums, and relationships involving multiplication, division, addition and subtraction.

The adjusted configuration satisfies all five specified structural properties. Under the randomisation model, none of 1{,}000{,}000 rearrangements of its own numerical values reached this score. Within the 64 permitted combinations of readings and adjustments, it was the only complete match; the row-wise product property alone selected this configuration. Both results are conditional on properties and adjustments identified during exploration of the artefact (Section~\ref{sec:stateva}).

Taken together, the documented patterns and the cross-column relationships developed through the reconstruction offer reasonable grounds for considering intentional arithmetic organisation a possibility. They motivate the hypothesis that the Ishango Bone may have served as a reference for demonstrating or teaching arithmetic relationships, potentially reflecting an early expression of mathematical reasoning. However, further material examination, experimental archaeological work and investigation of other marked artefacts are needed to assess the historical significance of this interpretation. The findings warrant such investigation and encourage closer examination of possible numerical organisation in other artefacts rather than its dismissal as incidental from the outset.

\section*{Data and Code Availability}
\label{subsec:DandCA}

\begin{sloppypar}
The Python scripts supporting the permutation analyses and the complete sensitivity analysis of readings and adjustments are archived on Zenodo at \url{https://doi.org/10.5281/zenodo.24023373}.

The maintained repository is publicly available at
\url{https://github.com/jbghia/ishango-bone-analysis}.

The analysis-code version used to obtain the reported results is identified by commit \texttt{d9ae017}: \url{https://github.com/jbghia/ishango-bone-analysis/commit/d9ae017eee18449f8dfda88cbb86f92b7fac23d6}.

The archived release includes subsequent README clarifications without changes to the analysis scripts.
\end{sloppypar}

The repository contains the reference configurations and property definitions used in the analyses, scripts for the sampled randomisation tests, and a complete enumeration of 64 specified combinations of readings and adjustments. The latter generates machine-readable CSV tables containing the configurations, individual and pairwise property outcomes, score distributions, and sensitivity to each reconstruction choice, together with a plain-text summary report. The repository README provides execution instructions and documents the scope and limitations of the analyses. All scripts use only the Python standard library.

\section*{AI Disclosure}

Claude (Anthropic, claude.ai) and Perplexity (Perplexity AI, perplexity.ai) assisted in developing the statistical analyses described in Section~\ref{sec:stateva}, including the composite scoring framework, the configuration-specific Monte Carlo randomisation analyses, the
complete sensitivity analysis of readings and adjustments, and the supporting Python scripts. The scripts are publicly available in the GitHub repository and archived on Zenodo, as detailed in the Data and Code Availability statement. The author reviewed and approved the methodology, code, reported results, and interpretation and takes responsibility for the final content.

\newpage

\printbibliography

\end{document}